\documentclass[12pt]{article}
\usepackage[authoryear,round]{natbib}
\usepackage{graphicx,bbm}
\usepackage{latexsym}
\usepackage[left=1in,top=1in,right=1in,bottom=1in]{geometry}\usepackage{amsmath,amssymb}
\usepackage{booktabs}
\usepackage{color}
\usepackage{tikz}
\definecolor{darkblue}{rgb}{0,0.4,0.9}
\definecolor{gray10}{rgb}{0.1,0.1,0.1}
\definecolor{gray20}{rgb}{0.2,0.2,0.2}
\definecolor{gray30}{rgb}{0.3,0.3,0.3}
\definecolor{gray40}{rgb}{0.4,0.4,0.4}
\definecolor{gray60}{rgb}{0.6,0.6,0.6}
\definecolor{gray80}{rgb}{0.8,0.8,0.8}
\definecolor{gray90}{rgb}{0.9,0.9,.9}
\definecolor{gray95}{rgb}{0.95,0.95,.95}
\definecolor{gray96}{rgb}{0.96,0.96,.96}
\definecolor{lgreen} {RGB}{180,210,100}
\definecolor{dblue}  {RGB}{20,66,129}
\definecolor{ddblue} {RGB}{11,36,69}
\definecolor{lred}   {RGB}{220,0,0}
\definecolor{nred}   {RGB}{224,0,0}
\definecolor{norange}{RGB}{230,120,20}
\definecolor{nyellow}{RGB}{255,221,0}
\definecolor{ngreen} {RGB}{98,158,31}
\definecolor{dgreen} {RGB}{78,138,21}
\definecolor{nblue}  {RGB}{28,130,185}
\definecolor{jblue}  {RGB}{20,50,100}
\definecolor{nnyellow}{RGB}{235,200,0}
\definecolor{purple}{RGB}{150, 0, 120}
\definecolor{sgGreen} {RGB}{20, 180, 50}
\definecolor{revised}{rgb}{0,0,0.9}

\newtheorem{theorem}{Theorem}
\newtheorem{corollary}{Corollary}
\newtheorem{lemma}{Lemma}

\newcommand{\nl}{\newline}
\newcommand{\pl}{\parallel}
\newcommand{\openr}{\hbox{${\rm I\kern-.2em R}$}}
\newcommand{\openn}{\hbox{${\rm I\kern-.2em N}$}}

\title{Efficient Estimation of Pathwise Differentiable Target Parameters with the Undersmoothed Highly Adaptive Lasso}
\author{Mark J. van der Laan$^1$\thanks{email: \texttt{laan@berkeley.edu}}, David Benkeser$^2$ and Weixin Cai$^1$ \\
\small{$^1$Division of Biostatistics, University of California, Berkeley} \\
\small{$^2$ Department of Biostatistics and Bioinformatics, Emory University}}

\begin{document}
\maketitle

\begin{abstract} 
We consider estimation of a functional  parameter of a 
realistically modeled data distribution based on observing independent and identically distributed observations. 
We define an $m$-th order Spline Highly Adaptive Lasso Minimum Loss Estimator (Spline HAL-MLE) of a functional parameter that is defined by minimizing the empirical risk function over an $m$-th order smoothness class of functions. We show that this $m$-th order smoothness class consists of all functions that can be represented as an infinitesimal linear combination of tensor products of $\leq m$-th order spline-basis functions, and involves assuming $m$-derivatives in each coordinate. By selecting $m$ with cross-validation we obtain a Spline-HAL-MLE that is able to adapt to the underlying unknown smoothness of the true function, while guaranteeing a rate of convergence faster than $n^{-1/4}$, as long as the true function is cadlag (right-continuous with left-hand limits) and has finite sectional variation norm. The $m=0$-smoothness class consists of all cadlag functions with finite sectional variation norm and corresponds with the original HAL-MLE defined in \citet{vanderLaan15}.

In this article we establish that this Spline-HAL-MLE yields an asymptotically efficient estimator of any smooth feature of the functional parameter under an easily verifiable global undersmoothing condition.  A sufficient condition for the latter condition is  that the minimum of the empirical mean of the selected basis functions is smaller than a constant times $n^{-1/2}$, which is not parameter specific. Therefore, the undersmoothing condition enforces the selection of the $L_1$-norm in the lasso to be large enough so that the fit includes  sparsely supported basis functions. We demonstrate our general result for the $m=0$-HAL-MLE of  the average treatment effect  and of the integral of the square of the data density. We also present simulations for these two examples confirming the theory.

\end{abstract}

\section{Introduction}
We consider the estimation problem in which we observe $n$ independent and identically distributed copies of a random variable with probability distribution known to be an element of an infinite-dimensional statistical model, while the goal is to estimate a particular smooth functional of the data distribution. It is assumed that the target parameter is a pathwise differentiable functional of the data distribution so that its derivative is characterized by its so called canonical gradient.

A regular asymptotically linear estimator is asymptotically efficient if and only if it is asymptotically linear with influence curve the canonical gradient  \citep{Bickeletal97} and a number of general methods for efficient estimation have been developed in the literature. If the model is not too large, then a regularized or sieve maximum likelihood estimator or minimum loss estimator (MLE) generally results in an efficient substitution estimator \citep{Newey94,vanderLaan06,vanderVaart98}. For a general theory on sieve estimation that also demonstrates sieve-based maximum likelihood estimators that are asymptotically efficient in large models, we refer to \citet{shen97,shen07}. These results generally require a sieve-based MLE that overfits the data (or equivalently, undersmooths the estimated functional parameter) and are only applicable for certain type of sieves. 

An alternative to undersmoothing is to use targeted estimator based on the canonical gradient, such as: the one-step estimator, which adds to an initial plug-in estimator the empirical mean of the canonical gradient at the estimated data distribution \citep{Bickeletal97};  an estimating equations-based estimator, which defines the estimator of the target parameter as the solution of an estimating equation with the estimated canonical gradient as estimating function \citep{Robins:Rotnitzky92,vanderLaan:Robins03}; and targeted minimum loss-estimation, which updates an initial estimator of the data distribution with an MLE of a least favorable parametric submodel through the initial estimator \citep{vanderLaan&Rubin06,vanderLaan08, vanderLaan&Rose11,vanderLaan&Gruber15}. By using an initial estimator of the relevant parts of the data distribution that converges w.r.t. $L^2$-type norm to the truth at a rate faster than $n^{-1/4}$, such as achieved with the HAL-MLE \citep{vanderLaan15,Benkeser&vanderLaan16}, in great generality, these three general procedures will result in an efficient estimator.
 
In this article we focus on a particular sieve MLE, which we call the HAL-MLE. The HAL-MLE is defined as the minimizer of an empirical mean of the loss function (e.g, log-likelihood loss) over a class of functions that can be arbitrarily well approximated by linear combinations of tensor products of univariate spline-basis functions, but where the $L_1$-norm of the coefficient vector is constrained. The target parameter is defined as a particular smooth real- or Euclidean-valued function of the functional parameter estimated by HAL-MLE, so that the HAL-MLE results in a plug-in estimator of the target parameter. In this case the sieve is indexed by a bound on the $L_1$-norm. By increasing this bound up to a large, finite value, the sieve approximates the total parameter space for the true functional parameter. If the goal is to estimate the functional itself, then the constraint on the $L_1$-norm is optimally chosen with cross-validation. In particular, the HAL-MLE described in \citet{vanderLaan15,Benkeser&vanderLaan16} selects the tuning parameter that minimizes the empirical mean of the loss function over the class of cadlag functions with finite sectional variation norm which can be approximated by an infinite linear combination of tensor product of indicator basis functions (i.e., a 0-th order spline-basis). In this case the $L_1$-norm of the coefficients equals the sectional variation norm of the function \citep{gill1995annalspoincare,vanderLaan15}.

The contributions of this article are two-fold. First, we generalize the 0-th order HAL-MLE to an $m$-th order Spline HAL-MLE in the class of $m$-times differentiable functions that can be approximated as a linear combination of tensor product of $\leq m$-th order spline basis functions with a finite $L_1$-norm of the coefficient vector. In this case, we refer to the $L_1$-norm of the coefficients as an $m$-th order sectional variation norm.  The algorithms for implementing these  $m$-th order Spline HAL-MLEs are identical across $m$ (just different basis functions) and can be based on implementations of the Lasso in the machine learning literature. One can now select both the bound on the $L_1$-norm and the smoothness degree $m$ with cross-validation, resulting in an estimator we call the smoothness-adaptive spline-HAL-MLE of the functional parameter. 

Second, we investigate whether and how an appropriately undersmoothed $m$-th order Spline HAL-MLE can be used to produce an efficient plug-in estimator of smooth functions of the functional parameter. There are essentially three key ingredients to establishing efficiency of a plug-in estimator: negligibility of the empirical mean of the canonical gradient, control of the second-order remainder, and asymptotic equicontinuity. For the first point, we argue that since the canonical gradient is a score, we essentially require that HAL-MLE solves a particular score equation. Because HAL-MLE is an MLE, it solves a large class of score equations, and we investigate whether these score equations might also approximate the particular score equation implied by the canonical gradient. In particular, we find that the larger the $L_1$-norm of the HAL-MLE, the more such score equations are generated and solved by the HAL-MLE. Therefore, one expects that by increasing the $L_1$-norm of the HAL-MLE, the linear span of equations solved by the HAL-MLE will approximate in first-order the canonical gradient score equation. However, another crucial condition for efficiency of a plug-in estimator is that a second-order remainder is $o_P(n^{-1/2})$, and we want to preserve the $n^{-1/4}$-rate of convergence of achieved by the HAL-MLE when the $L_1$-norm is selected with cross-validation. Fortunately, the  rate of the HAL-MLE is not affected by the size of the $L_1$-norm as long as it  remains bounded and, for $n$ large enough, exceeds the $m$-th order sectional variation norm of the true function. Similarly, the asymptotic equicontinuity condition for efficiency of a plug-in estimator will also be satisfied for any bounded $L_1$-norm, since the class of cadlag functions with a finite $0$-order sectional variation norm is a Donsker class.  In fact, one can prove that this $L_1$-norm is allowed to slowly converge to infinity as sample size increases without affecting the asymptotic equicontinuity  condition and the $n^{-1/4}$-rate of convergence of the HAL-MLE. Taken together, our analysis highlights that when selected the level of undersmoothing of a HAL-MLE, one wants to undersmooth enough to solve the efficient score equation up to an appropriate level of approximation, but in order to reasonable finite-sample performance one should not undersmooth beyond that level. 

This discussion highlights the need to establish empirical criterion by which the level of undersmoothing may be chosen to appropriately satisfy the conditions required of an efficient plug-in estimator. In particular, we an easily verifiable global undersmoothing condition, which is satisfied for example, if the minimum of the empirical mean of the basis functions that receive non-zero coefficient is smaller than a constant times $n^{-1/2}$. This condition essentially enforces the selection of the $L_1$-norm in the Lasso to be large enough so that the fit includes  sparsely supported basis functions. We also discuss alternative practical criterion for selecting the level of undersmoothing. We demonstrate our result in practice for the $m=0$-HAL-MLE of the average treatment effect in a nonparametric model, and for estimation of the integral of the square of the data density. 


This article is organized as follows. 
In the next Section \ref{section2} we define the $m$-th order HAL-MLE; a formal proof of the representation theorem is provided in the Appendix. 
In Section \ref{section3} we establish our main theorem providing the undersmoothing conditions under which the $m$-th order Spline-HAL-MLE is asymptotically efficient for any pathwise differentiable parameter.
In Section \ref{section4} we apply our theorem to the ATE example providing a theorem for this particular nonparametric estimation problem. In Section \ref{section5} we apply our theorem to a nonparametric estimation problem with target parameter the integral of the square of the data density.  
In Section \ref{section6} we demonstrate a simulation study for both examples, providing a practical verification of our theoretical results. 

\section{Defining functional estimation problem, and Spline-HAL-MLE}\label{section2}

\subsection{Functional estimation problem}
Suppose we observe $O_1,\ldots,O_n\sim_{iid} P_0\in {\cal M}$, where $O$ is a Euclidean random variable of dimension $k_1$ with support contained in $[0,\tau_o]\subset\openr^{k_1}$.
Let $Q:{\cal M}\rightarrow Q({\cal M})=\{Q(P):P\in {\cal M}\}$ be a functional parameter. It is assumed that there exists a loss function $L(Q)$  so that $P_0L(Q(P_0))=\min_{P\in {\cal M}}P_0L(Q(P))$, where we use the notation $Pf\equiv \int f(o)dP(o)$. Thus, $Q(P)$ can be defined as the minimizer of the risk function $Q\rightarrow PL(Q)$ over all $Q$ in the parameter space. Let $d_0(Q,Q_0)\equiv P_0L(Q)-P_0L(Q_0)$ be the loss-based dissimilarity.
We assume that $M_{20}\equiv \sup_{P\in {\cal M} } P_0\{L(Q(P))-L(Q_0)\}^2/d_0(Q(P),Q_0)<\infty$, and $M_1\equiv \sup_{0,P\in {\cal M}}\mid L(Q(P))(o)\mid<\infty$, thereby guaranteeing good behavior of the cross-validation selector \citep{vanderLaan&Dudoit03,vanderVaart&Dudoit&vanderLaan06,vanderLaan&Dudoit&vanderVaart06,vanderLaan&Polley&Hubbard07,Chpt3}.

{\bf Parameter space for functional parameter $Q$: Cadlag and  uniform bound on sectional variation norm.}
We assume that the parameter space $Q({\cal M})$ is a collection  of multivariate real valued cadlag functions on a cube $[0,\tau]\subset\openr^k$
with finite sectional variation norm $\pl Q(P)\pl_v^*< C^u$ for some $C^u<\infty$. That is, for all $P$, $Q(P)$ is a $k$-variate real valued cadlag function on $[0,\tau]\subset\openr^k_{\geq 0}$ with $\pl Q(P)\pl_v^*<C^u$, where the sectional variation norm is defined by \[
\pl Q\pl_v^*\equiv Q(0)+\sum_{s\subset\{1,\ldots,k\}}\int_{[0_s,\tau_s]}\mid dQ_s(u_s)\mid.\]
For a given subset $s\subset\{1,\ldots,k\}$,  $Q_s:(0_s,\tau_s]\rightarrow\openr$ is defined by $Q_s(x_s)=Q(x_s,0_{-s})$. That is, $Q_s$ is  the $s$-specific section of $Q$ which sets the coordinates in the complement of  subset $s\subset\{1,\ldots,k\}$ equal to 0. Since $Q_s$ is right-continuous with left-hand limits and has a finite variation norm over $(0_s,\tau_s]$, it generates a finite measure, so that the integrals w.r.t. $Q_s$ are indeed well defined. For a given vector $x\in[0,\tau]$, we define
$x_s=(x(j):j\in s)$. Sometimes, we will also use the notation $x(s)$ for $x_s$. 

Note also that $[0,\tau]=\{0\}\cup\left(  \cup_s (0_s,\tau_s]\right)$ is partitioned in the singleton $\{0\}$, the $s$-specific left-edges $(0_s,\tau_s]\times \{0_{-s}\}$ of cube $[0,\tau]$, and, in particular, the full-dimensional inner set $(0,\tau]$ (corresponding with $s=\{1,\ldots,k\}$). Therefore, the above sectional variation norm equals the sum over all subsets $s$ of the variation norm of the $s$-specific section over its $s$-specific edge. 
It is also important to note that any cadlag function $Q$ with finite sectional variation norm can be represented as
\[
Q(x)=Q(0)+\sum_{s\subset \{1,\ldots,k\}}\int_{(0_s,x_s]} dQ_s(u_s).\] 
That is, $Q(x)$ is a sum of integrals up till $x_s$ over all the $s$-specific edges w.r.t.  the measure generated by the corresponding $s$-specific section $Q_s$.  
We will refer to $Q_s$ as a cadlag function as well as a measure.
We note that this representation represents $Q$ as an infinitesimal linear combination of indicator basis functions $x\rightarrow \phi_{s,u_s}(x)\equiv I(x_s\geq u_s)$ indexed by knot-point $u_s$ with coefficient $dQ_s(u_s)$:
\[
Q(x)=Q(0)+\sum_{s\subset\{1,\ldots,k\}}\int \phi_{s,u_s}(x) dQ_s(u_s).\]
Note that the $L_1$-norm of the coefficients in this representation is precisely the sectional variation norm $\pl Q\pl_v^*$.

\subsection{$m$-th order spline smoothness class and its $m$-th order spline representation}

{\bf Iterative definition of relevant $m$-th order derivatives of $Q$:}
Our $m$-th order smoothness class relies on the existence of certain $m$-th order derivatives. We will  now define these $m$-th order derivatives through recursion.
For a function $Q$, we define the $s$-specific section $Q_s(x_s)=Q(x_s,0_{-s})$.
The first order derivative $Q_s^1$ of $Q_s$ is defined as a density of $Q_s$ w.r.t. Lebesgue measure so that $dQ_s(u_s)=Q_s^1(u_s)du_s$.
Given the set of first order derivatives $\{Q_s^1:s\subset\{1,\ldots,k\}\}$ indexed by all subsets $s\subset\{1,\ldots,k\}$, we will now define the set of second order derivatives
$\{Q_{s,s_1}^2: s,s_1\subset s\}$, indexed by all subsets $s\subset\{1,\ldots,k\}$ and all its subsets $s_1\subset s$. Given the function $Q_s^1$ and $s_1\subset s$, we define its $s_1$-specific section $Q_{s,s_1}^1(x_{s_1})=Q_s^1(x_{s_1},0_{-s_1})$. 
The second order derivative $Q_{s,s_1}^2$ is defined  as the density of $Q_{s,s_1}^1$ w.r.t. Lebesgue measure so that  $Q_{s,s_1}^1(du_{s_1})=
Q_{s,s_1}^2(u_{s_1})du_{s_1}$.
This defines now $\{Q_{\bar{s}(1)}^2:\bar{s}(1)\}$, where $\bar{s}(1)=(s,s_1)$ and it varies over all $s\subset\{1,\ldots,k\}$ and subsets $s_1\subset s$.

Let $m=2$. Given the set of $m$-th order derivatives $\{Q_{\bar{s}(m-1)}^m:\bar{s}(m-1)\}$, we will now define 
$\{\bar{Q}_{\bar{s}(m)}^{m+1}:\bar{s}(m)\}$. We are reminded again that $\bar{s}(m)=(s,s_1,\ldots,s_m)$ is a sequence of nested subsets, $s_m\subset s_{m-1}\subset\ldots\subset s$. We also note that $Q_{\bar{s}(m-1)}^m(x_{s_{m-1}} )$ is only a function of coordinates in $s_{m-1}$, since all other coordinates have been set to zero through the earlier sections implied by $s_{m-2},\ldots,s$.
Given the function $Q_{\bar{s}(m-1)}^{m}$, we define its $s_m$-specific section $Q_{\bar{s}(m)}^m(x_{s_m})=Q_{\bar{s}(m-1)}^m(x_{s_{m}},0_{-s_{m}})$ that sets the coordinates in $s_{m+1}/s_m$ equal to zero.  The $m+1$-th order derivative $Q_{\bar{s}(m)}^{m+1}$ is defined  as the density of $Q_{\bar{s}(m)}^m$ w.r.t. Lebesgue measure so that  $Q_{\bar{s}(m)}^m(du_{s_m})=
Q_{\bar{s}(m)}^{m+1}(u_{s_m})du_{s_m}$.

{\bf $m$-th order sectional variation norm:}
The $m$-th order sectional variation norm is defined as:
\[
\begin{array}{l}
\pl Q\pl_v^{*,m}\equiv \mid Q(0)\mid+\sum_{j=1}^{m-1}\sum_{\bar{s}(j)} \mid Q_{\bar{s}(j)}^{j+1}(0_{s_{j}})\mid 
+\sum_{\bar{s}(m)}\int_{(0_{s_m},\tau_{s_m}]}\mid dQ^m_{\bar{s}(m)}(z_{s_m})\mid .
\end{array}
\]

{\bf Iterative definition of relevant $m$-th order spline basis functions:}
We first describe an iterative procedure that allows us to define the relevant $m$-th order spline basis functions whose linear span generates the $m$-th order smoothness class defined below.

For  $s\subset\{1,\ldots,k\}$, we define the $0$-order spline basis functions $\phi_{s,x_s}(u_s)=I(x_s\geq u_s)$ indexed by knot point $u_s$.
For $s_1\subset s$, we define the first order spline  basis functions
$\phi_{s,s_1,x_{s}}(u_{s_1})=\prod_{j\in s_1}(x(j)-u(j))I(u(j)\leq x(j))\prod _{j\in s/s_1}x_j$, indexed by knot-point  $u_{s_1}$. We also define \[
\phi_{s,\emptyset,x_s}(0_{s_1})=\prod_{j\in s}x_j,\]
which corresponds with setting $s_1$ equal to empty set in definition of $\phi_{s,s_1,x_s}(u_{s_1})$.

For a given $\bar{s}(m)$ and corresponding basis function $\phi_{\bar{s}(m),x_s}(u_s)=\prod_{j\in s}\phi_{j,\bar{s}(m),x_j}(u_j)$, we note that it is a tensor product of univariate basis functions $\phi_{j,\bar{s}(m),x_j}(u_j)$ over the components $j\in s$. 
Let $m=1$.
Given $\bar{s}(m)$, and given $s_{m+1}\subset s_m$, we define the $m+1$-th order spline basis functions as
\begin{eqnarray*}
\phi_{\bar{s}(m+1),x_s}(z_{s_{m+1}})
&\equiv& \prod_{j\in s_{m+1}}\int_{(z_j,x_j]}\phi_{j,\bar{s}(m),x_j}(y_j)dy_j\prod_{j\in s_{m}/s_{m+1}}\int_{(0,x_j]}\phi_{j,\bar{s}(m),x_j}(y_j)dy_j\\
&&
\prod_{j\in s/s_m}\phi_{j,\bar{s}(m),x_j}(0).
\end{eqnarray*}
We also define \[
\phi_{\bar{s}(m),\emptyset,x_s}(0_{s_{m+1}})=\prod_{j\in s_{m}}\int_{(0,x_j]}\phi_{j,\bar{s}(m),x_j}(y_j)dy_j \prod_{j\in s/s_m}\phi_{j,\bar{s}(m),x_j}(0)\]
 by setting $s_{m+1}$ equal to empty set, and knot point $z_{s_{m+1}}=0_{s_{m+1}}$ in the definition of $\phi_{\bar{s}(m+1),x_s}(z_{s_{m+1}})$.
Note that for each $j\in s_{m+1}$, the previous $m$-th order basis function is smoothed by integrating it from a knot point $z_j$ till $x_j$; for $j\in s_m/s_{m+1}$ the previous $m$-th order basis function is smoothed by integrating it  from $0$ till $x_j$; and, finally, for $j\in s/s_m$, the $m$-th order basis function is untouched.

{\bf $m$-th order spline smoothness class:}
Let $D^m[0,\tau]$ be the space  of cadlag functions  $f:[0,\tau]\rightarrow\openr$ for which the $m$-th order derivatives $\{f_{\bar{s}(m)}^m:\bar{s}(m)\}$ exist, and $m$-th order sectional variation norm is bounded, $m=0,1,\ldots$. Let $D^m_C[0,\tau]=\{f\in D^m[0,\tau]: \pl f\pl_v^{*,m}<C\}$ be the subset that enforces the $m$-th order sectional variation norm to be bounded by $C$. 
We refer to $D^m[0,\tau]$ as the $m$-th order spline  smoothness class. 

In the Appendix we establish the following representation for this smoothness class $D^m[0,\tau]$.
For any function $Q\in D^m[0,\tau]$ (i.e., finite $m$-th order sectional variation norm), we have
\begin{eqnarray*}
Q(x)&=&Q(0)+\sum_{j=0}^{m-1}\sum_{\bar{s}(j)}Q_{\bar{s}(j)}^{j+1}(0_{s_j})\phi_{\bar{s}(j),\emptyset,x_s}(0)\\
&&+\sum_{\bar{s}(m)}\int \phi_{\bar{s}(m),x_s}(z_{s_m})dQ^m_{\bar{s}(m)}(z_{s_m}).
\end{eqnarray*}
Just as for the $m=0$-smoothness class of cadlag functions with finite sectional variation norm,   a $Q\in D^m[0,\tau]$ is represented by  an infinitesimal linear combination of a collection of basis functions  given by
\[
{\cal F}_{b}^m\equiv \cup_{j=1}^{m-1}\{x\rightarrow \phi_{\bar{s}(j),\emptyset,x_s}(0):\bar{s}(j)\}\cup\{
x\rightarrow \phi_{\bar{s}(m),x_s}(z_{s_m}): \bar{s}(m),z_{s_m}\in (0_{s_m},\tau_{s_m}]\}.
\]
Let ${\cal S}_j$ be the collection of nested  $\bar{s}(j)$ of subsets with $s_j\subset s_{j-1}\ldots\subset s$, and $s\subset\{1,\ldots,k\}$. Define
\begin{eqnarray*}
{\cal J}^m_d&=&\cup_{j=1}^{m-1}{\cal S}_j\\
{\cal J}^m_c&=&\{\bar{s}(m),z_{s_m}: \bar{s}(m)\in {\cal S}_m,z_{s_m}\in (0_{s_m},\tau_{s_m}]\}\\
{\cal J}^m&=&{\cal J}^m_d\cup{\cal J}^m_c.
\end{eqnarray*}
Notice that ${\cal J}^m$ is a union of a finite index set ${\cal J}^m_d$ and an infinite set ${\cal J}^m_c$.
In addition, for notational, convenience, for each index $i\in {\cal J}^m$, let $\phi_i$ denote  the corresponding basis function as defined above (i.e., either of type $\phi_{\bar{s}(j),\emptyset,x_s}(0)$ or $\phi_{\bar{s}(m),x_s}(z_{s_m})$). 
Then, we can represent the class of basis functions as follows:
\[
{\cal F}_b^m=\left\{\phi_{i}:i\in {\cal J}^m_d\}\cup\{\phi_i:i\in {\cal J}^m_c\right\}.\]
In addition, we have 
\[
Q(x)=\sum_{i\in {\cal J}^m_d}\phi_i(x)\beta(i)+\int_{i\in {\cal J}^m_c}\int_i \phi_i(x) d\beta(i).\]
For $i\in {\cal J}^m_c$ with $i=(\bar{s}(m),z_{s_m})$ we define $d\beta(i)=dQ^m_{\bar{s}(m)}(z_{s_m})$, and for $i\in {\cal I}^m_d$ with $i=\bar{s}(j)$, we have $\beta(i)=Q_{\bar{s}(j)}^{j+1}(0_{s_j})$.

We could now represent 
\[
D^m[0,\tau]=\left\{\sum_{i\in {\cal J}^m_d}\beta(i)\phi_i+\int_{i\in {\cal J}^m_c} \phi_i d\beta(i):\pl \beta\pl_{1}<\infty\right\},\]
while
\[
D^m_C[0,\tau]=\left\{\sum_{i\in {\cal J}^m_d}\beta(i)\phi_i+\int_{i\in {\cal J}^m_c} \phi_i d\beta(i):\pl \beta\pl_{1}<C\right\},\]
where 
\[
\pl \beta\pl_1=\sum_{i\in {\cal J}^m_d}\mid \beta(i)+\int_{i\in {\cal J}^m_c} \mid d\beta(i)\mid .\]
{\bf Another convenient notation for spline representation:}
In order to emphasize that each of the basis functions $\phi_{\bar{s}(m),x_s}(u_{s_m})$ and
$\phi_{\bar{s}(j),\emptyset,x_s}(0)$ is a tensor product of $\leq m$-th order univariate spline-basis functions over components $j\in s$ ($s$ being first component of $\bar{s}(m)$), we will also use the notation $\phi_{s,j}$, where  $s\subset\{1,\ldots,k\}$ and $j$ represents an index in ${\cal J}^m$ with first subset equal to $s$. Let ${\cal J}^m(s)={\cal J}^m_d(s)\cup{\cal J}^m_c(s)$ are all elements in the index set ${\cal J}^m$ for which the first subset in its subset-vector equals $s$.
With this notation, we can write
\begin{equation}\label{reprparspace}
D^m[0,\tau]=\left\{\sum_{s,i\in {\cal J}^m_d(s)}\beta(s,i)\phi_{s,i}+\sum_s\int_{i\in {\cal J}^m_c(s)}\phi_{s,i}\beta(s,di):\pl \beta\pl_1<\infty\right\},\end{equation}
and similarly for $D^m_C[0,\tau]$.

\subsection{Definition of $m$-th order Spline-HAL-MLE}
Recall ${\cal Q}(C^u)=\{Q\in D[0,\tau]: \pl Q\pl_v^*<C^u\}$ be the class of cadlag functions which  with sectional variation norm bounded by $C^u$. 
Let $C_{0m}\equiv \pl Q_0\pl_v^{*,m}$ be the $m$-th order sectional variation norm of $Q_0$, and let $C^u_m$ be an upper bound guaranteeing that $C_{0m}<C^u_m$.
For a constant $C<C^u_m$,  consider the $m$-th order spline-function class ${\cal Q}^m(C)\equiv \{Q\in D^m[0,\tau]:\pl Q\pl_v^{*,m}<C,\pl Q\pl_v^*<C^u\}\subset {\cal Q}(C^u)$.  For a data adaptive selector $C_n$, we define
 \begin{equation}\label{unrestrictedHAL}
 Q_n^m\equiv \arg\min_{Q\in {\cal Q}^m(C_n)}P_n L(Q)\end{equation}
 be the $m$-th order Spline HAL-MLE. We will restrict the minimization to $Q$ for which, for all subset vectors $\bar{s}(m)$,  $dQ^m_{\bar{s}(m)}$ is a discrete measure with a finite support
 $\{z_{\bar{s}(m),j}:j=1,\ldots,n_{\bar{s}(m)}\}$. That is,  for each $\bar{s}(m)$,  the $m$-th order derivative $Q^m_{\bar{s}(m)}$  of $Q$ is absolutely continuous w.r.t. a discrete counting measure $\mu_{n,\bar{s}(m)}$. We will denote this form of absolute continuity with 
 $Q^m\ll^* \mu_n$. 
 Thus, $m$-th order HAL-MLE then becomes
 \[
 Q_n^m\equiv \arg\min_{Q\in {\cal Q}^m(C_n),Q^m\ll^*  \mu_n} P_n L(Q).\]
 Our $\leq m$-th  order spline representation for functions in $D^m[0,\tau]$ shows that all $Q$ with $Q^m\ll^*\mu_n$ are represented by a finite dimensional linear combination of basis functions indexed by  ${\cal J}^m(\mu_n)={\cal J}^m_d\cup  {\cal J}^m_c(\mu_n)$ for a finite subset ${\cal J}^m_c(\mu_n)\subset {\cal J}^m_c$.  
 Therefore, in this case the $m$-th order Spline-HAL MLE can be represented as  $Q_n^m=\sum_{j\in {\cal J}^m(\mu_n)} \beta_n^m(j)\phi_j$, where \[
 \beta_n^m\equiv \arg\min_{\beta,\pl \beta\pl_1\leq C_n}L\left (\sum_{j\in {\cal J}^m(\mu_n)} \beta(j)\phi_j\right ).\]
 Note that $Q_n^m=\hat{Q}^m(P_n)$ is the realization of a  mapping from the empirical probability measure to the parameter space. 
  
  As noted earlier, the data adaptive selector $C_n$ might be selected larger or equal than cross-validation selector $C_{n,m}=\arg\min_C E_{B_n}P_{n,B_n}^1 L(\hat{Q}^m(P_{n,B_n}^0))$,
  where $B_n\in \{0,1\}^n$ represents a random sample split (e.g, $V$-fold cross-validation)  into a training sample $\{i:B_n(i)=0\}$ and validation sample $\{i:B_n(i)=1\}$, while $P_{n,B_n}^0$ and $P_{n,B_n}^1$ are the corresponding empirical probability measures.
 One wants that $C_n\geq C_{0m}$ for $n$ large enough, so that  $Q_0\in {\cal Q}^m(C_n)$.

 Typically, one is able to prove that the unrestricted MLE (\ref{unrestrictedHAL}) will be discrete on a support in which case our $\mu_n$-discretization does not restrict the definition of the HAL-MLE.
 Generally, if $O$ includes observing $X$ where $L(Q)(0)$ depends on $Q$ through $Q(X)$, we recommend to select the support  of $dQ^m_{\bar{s}(m)}$ as a subset (or whole set) of the observed data $X_i(s_m)$, $i=1,\ldots,n$.


\paragraph{Smoothness adaptive spline HAL-MLE}
Suppose now that we select $m_n$ with the cross-validation selector. In addition, assume that each smoothness class allows for a unique rate of convergence of the corresponding spline HAL-MLE w.r.t. loss-based dissimilarity. Due to the asymptotic equivalence of the cross-validation selector with the oracle selector, it then follows that $P(m_n=m_0)\rightarrow 1$ where $m_0$ is the unknown true maximal smoothness of the true $Q_0$. See Appendix B for a formal statement. As a consequence, this smoothness adaptive spline HAL-MLE achieves the rate of convergence of the $m_0$-th smoothness class (i.e., it is minimax adaptive). The asymptotic efficiency of a smoothness adaptive spline HAL-MLE (using undersmoothing of each $m$-th order Spline HAL-MLE) follows from the asymptotic efficiency of the $m_0$-th order Spline HAL-MLE as established in the next section. 


\section{Efficiency of $m$-th order Spline-HAL MLE  for pathwise differentiable  target parameters}\label{section3}
\subsection{Defining the efficient estimation problem and plug-in HAL-MLE}
Let $\Psi:{\cal M}\rightarrow\openr^d$ be the $d$-dimensional statistical target parameter of interest of the data distribution. 
We assume that $\Psi$ is pathwise differentiable at any $P\in {\cal M}$ with canonical gradient $D^*(P)$. For a pair $P,P_0\in {\cal M}$, the exact second order remainder is defined by \[
R_2(P,P_0)\equiv \Psi(P)-\Psi(P_0)+P_0D^*(P).\]  

{\bf Relevant functional parameter  and its loss function:}
Let $Q:{\cal M}\rightarrow Q({\cal M})=\{Q(P):P\in {\cal M}\}$ be a functional parameter such that $\Psi(P)=\Psi_1(Q(P))$ for some $\Psi_1$. It is assumed that $Q$ is a functional parameter with parameter space $Q({\cal M})\subset {\cal Q}(C^u)=D_{C^u}[0,\tau]$ as defined above in Section \ref{section2}. Note that the model ${\cal M}$ does not make any smoothness assumptions on $Q$ beyond that it is a cadlag function with sectional variation norm bounded by $C^u$.
In particular, we have an $m$-th order Spline-HAL-MLE with established rate of convergence for $d_0(Q_n^m,Q_0)=o_P(n^{-1/2})$ for $m\geq m_0$, where $m_0$ is the smoothness degree of $Q_0$.
In addition, due to the asymptotic equivalence of the cross-validation selector with the oracle selector, for the smoothness adaptive $Q_n^{m_n}$, where $m_n$ is the cross-validation selector of $m$,  we have $P(m_n=m_0)\rightarrow 1$ and $d_0(Q_n^{m_n},Q_0)=o_P(n^{-1/2})$.

{\bf Nuisance parameter for canonical gradient:}
 Let $G:{\cal M}\rightarrow{\cal G}$ be a functional nuisance parameter so that $D^*(P)$ only depends on $P$ through $(Q(P),G(P))$, and the remainder $R_2(P,P_0)$ only involves differences between $(Q,G)$ and $(Q_0,G_0)$:
 \[
 \mbox{ $D^*(P)=D^*(Q(P),G(P))$, while
 $R_2(P,P_0)=R_{20}((Q,G),(Q_0,G_0))$.}
 \]
 Here $R_{20}$ could have some remaining dependence on $P_0$ and $P$, and ${\cal G}=G({\cal M})$ is the parameter space for $G$.

{\bf Canonical gradient of target parameter in tangent space of loss function:}
We also assume that this loss function $L(Q)$ is such that there exists a class of submodels $\{Q_{\epsilon}^h: \epsilon\}\subset Q({\cal M})$ indexed by a choice $h\in {\cal H}^1$,  through $Q$ at $\epsilon=0$, so that for any $G\in {\cal G}$, one of these directions $h$ generates a score that equals the canonical gradient $D^*(Q,G)$ at $(Q,G$): \[
\left . \frac{d}{d\epsilon}L(Q_{\epsilon}^h)\right |_{\epsilon =0}=D^*(Q,G).
\]
Since the canonical gradient is an element of the tangent space and thereby typically a score of a submodel,  this generally holds for $Q$  defined as the density of $P$ and the log-likelihood loss $L(Q)=-\log Q$.  However, for any $Q$ there are  typically more direct loss functions $L(Q)$, so that the loss-based dissimilarity $d_0(Q,Q_0)=P_0L(Q)-P_0L(Q_0)$  directly measures a dissimilarity between $Q$ and $Q_0$, for which this condition holds as well.

{\bf $m$-th order HAL-MLE:}
Let $m\in \{0,1\ldots\}$ be fixed. In this section, we are concerned with analyzing the plug-in estimator $\Psi(Q_n^m)$ of $\Psi(Q_0)$, where $Q_n^m$ is the $C_n^m$-tuned $m$-th order Spline-HAL-MLE $\hat{Q}^m(P_n)=\hat{Q}_{C_n^m}^m(P_n)$, which minimizes the empirical risk over ${\cal Q}^m(C_n^m)$. We assume here that $m\geq m_0$, so that $Q_0\in {\cal Q}^m(C_m^u)$, even though that is not an assumption in our model ${\cal M}$. 
Recall our assumptions on the parameter space ${\cal Q}={\cal Q}(C^u)$.
 We assume that ${\cal Q}$ is defined such that $Q_n^m$ is in the interior of the model based parameter space ${\cal Q}$ (so that there are submodels through $Q_n^m$ that generate the tangent space and the canonical gradient), even though $Q_n^m$ is typically on the edge of the parameter subspace ${\cal Q}^m(C_n^m)\subset {\cal Q}=Q({\cal M})$ over which the estimator is minimizing the empirical risk.
It is understood that verification of our conditions might require using a $C_n^m$ different from the cross-validation selector.

\paragraph{Remark: Target parameter could be component of real target parameter.}
In many situations the real target parameter is a $P\rightarrow \Psi(Q_1(P),Q_2(P))$ for two (or more) functional parameters $Q_1$ and $Q_2$.
One could apply our efficiency theorem below to the target parameter $\Psi_{Q_{10}}(Q_2)=\Psi(Q_{10},Q_2)$ and $\Psi_{Q_{20}}(Q_1)=\Psi(Q_1,Q_{20})$ treating the indices $Q_{10}$ and $Q_{20}$ as known, and $m$-th order Spline-HAL-MLEs $Q_{1n}$ and $Q_{2n}$ of $Q_{10}$ and $Q_{20}$, respectively.
Application of our theorem to these two cases then proves that $\Psi(Q_{10},Q_{2n})$ and $\Psi(Q_{1n},Q_{20})$ are both asymptotically efficient, if both HAL-MLEs are appropriately tuned w.r.t. $m$-th order sectional variation norm bound. 
Since \[
\Psi(Q_{1n},Q_{2n})-\Psi(Q_{10},Q_{20})=\Psi(Q_{1n},Q_{2n})-\Psi(Q_{10},Q_{2n})+\Psi(Q_{10},Q_{2n})-\Psi(Q_{10},Q_{20}),\]
 this then also establishes asymptotic efficiency of $\Psi(Q_{1n},Q_{2n})$ as estimator of $\Psi(Q_{10},Q_{20})$, under the condition that
\[
\Psi(Q_{1n},Q_{2n})-\Psi(Q_{10},Q_{2n})-(\Psi(Q_{1n},Q_{20})-\Psi(Q_{10},Q_{20}))=o_P(n^{-1/2}).
\]
This latter term can be viewed as a second order difference of $(Q_{1n},Q_{2n})$ and $(Q_{10},Q_{20})$ so that the latter condition will generally hold by using the already established  rates of convergence $o_P(n^{-1/2})$ w.r.t. risk based dissimilarity  for $Q_{1n}$ and $Q_{2n}$. The above immediately generalizes to the case that the target parameter is a function of more than two  $Q$-components.

\subsection{The Spline-HAL MLE solves efficient influence curve equation by including sparse basis functions}

Let $Q_n=\arg\min_{Q\in {\cal Q}^m(C_n)}P_n L(Q)$ be the $m$-th order Spline HAL-MLE, suppressing the  dependence of $Q_n$ and $C_n$ on $m$. 
The following theorem establishes that, if $m\geq m_0$, then  $\Psi(Q_n)$ is asymptotically efficient for $\Psi(Q_0)$ for large enough $C_n$, and some weak conditions specific towards the target parameter.
It relies on the following definitions that also provide the basis of the proof of the theorem.
\nl
{\bf Definitions:} \ \nl
{\bf $\bullet$}
Recall we can represent $Q_n=\arg\min_{Q\in {\cal Q}^m(C_n)}P_n L(Q)$ as follows:
\[
Q_n(x)=I(Q_n)(x)+\sum_{\bar{s}(m)}\int_{(0_{s_m},x_{s_m}] } \phi_{\bar{s}(m),x_s}(u_{s_m})dQ^m_{n,\bar{s}(m)}(u_{s_m}),\]
where
\[
\begin{array}{l}
I(Q_n)(x)=Q_n(0)+\sum_{j=0}^{m-1}\sum_{\bar{s}(j)}Q_{n,\bar{s}(j)}^{j+1}(0_{s_j})\phi_{\bar{s}(j),\emptyset,x_s}(0_s).
\end{array}
\]
{\bf $\bullet$} Consider the family of paths $\{Q_{n,\epsilon}^h:\epsilon\in (-\delta,\delta)\}$ through $Q_n$ at $\epsilon =0$ for arbitrarily small $\delta>0$, indexed by any uniformly bounded $h$, defined by
\begin{equation}\label{familypaths}
Q_{n,\epsilon}^h(x)=I(Q_{n,\epsilon}^h)(x)+\sum_{\bar{s}(m)}\int_{(0_{s_m},x_{s_m}] } \phi_{\bar{s}(m),x_s}(u_{s_m})(1+\epsilon h(\bar{s}(m),u_{s_m}))dQ^m_{n,\bar{s}(m)}(u_{s_m}),\end{equation}
where
\begin{eqnarray*}
I(Q_{n,\epsilon}^h)(x)&=&(1+\epsilon h(0))Q_n(0)
+\sum_{j=0}^{m-1}\sum_{\bar{s}(j)}\phi_{\bar{s}(j),\emptyset,x_s}(0_s)
(1+\epsilon h(\bar{s}(j),0_{s_j}))Q_{n,\bar{s}(j)}^{j+1}(0_{s_j}).
\end{eqnarray*}
{\bf $\bullet$} Let 
\begin{eqnarray*}
r(h,Q_n)&\equiv& I(h,Q_n)+\sum_{\bar{s}(m)}\int_{(0_{s_m},\tau_{s_m}] }h(\bar{s}(m),u_{s_m})\mid dQ^m_{n,\bar{s}(m)}(u_{s_m})\mid,\end{eqnarray*}
where
\[
\begin{array}{l}
I(h,Q_n)=h(0))\mid Q_n(0)\mid +\sum_{j=0}^{m-1}\sum_{\bar{s}(j)} h(\bar{s}(j),0_{s_j})  \mid Q_{n,\bar{s}(j)}^{j+1}(0_{s_j})\mid .
\end{array}
\]
{\bf $\bullet$} For any uniformly bounded $h$ with $r(h,Q_n)=0$ we have that for a small enough $\delta>0$ $\{Q_{n,\epsilon}^h:\epsilon\in (-\delta,\delta)\}\subset {\cal Q}^m(C_n)$.\nl
{\bf $\bullet$} Let $S_h(Q_n)=\left . \frac{d}{d\epsilon}L(Q_{n,\epsilon}^h)\right |_{\epsilon=0}$ be the score of this $h$-specific submodel.
\nl
{\bf $\bullet$}
Consider the set of scores 
\begin{equation}\label{scores}
{\cal S}(Q_n)=\{S_h(Q_n)=\frac{d}{dQ_n}L(Q_n)(f(h,Q_n)):\pl h\pl_{\infty}<\infty\},
\end{equation} where
\begin{eqnarray*}
f(h,Q_n)(x)&\equiv &\left . \frac{d}{d\epsilon}Q_{n,\epsilon}^h\right |_{\epsilon =0}(x)\\
&=&f_1(h,Q_n)(x)+\sum_{\bar{s}(m)}\int_{(0_{s_m},x_{s_m}] } \phi_{\bar{s}(m),x_s}(u_{s_m})h(\bar{s}(m),u_{s_m})dQ^m_{n,\bar{s}(m)}(u_{s_m})\\
f_1(h,Q_n)(x)&=&\left . \frac{d}{d\epsilon}I(Q_{n,\epsilon}^h)(x)\right |_{\epsilon =0}\\
&=&h(0))Q_n(0)+\sum_{j=0}^{m-1}\sum_{\bar{s}(j)}\phi_{\bar{s}(j),\emptyset,x_s}(0_s)h(\bar{s}(j),0_{s_j})Q_{n,\bar{s}(j)}^{j+1}(0_{s_j}),
\end{eqnarray*}
This is the set of scores generated by the above class of paths if we do not enforce constraint $r(h,Q_n)=0$.
\nl
{\bf $\bullet$}  
We have that $Q_n$ solves the score equations $P_n S_h(Q_n)=0$ for any uniformly bounded $h$ satisfying $r(h,Q_n)=0$.
\nl {\bf $\bullet$} Let $D^*_n(Q_n,G_0)\in {\cal S}(Q_n)$ be an approximation of $D^*(Q_n,G_0)$ that is contained in this set of scores ${\cal S}(Q_n)$.
\nl {\bf $\bullet$} We also consider a special case in which $D^*_n(Q_n,G_0)=D^*(Q_n,G_{0n})$ for an approximation $G_{0n}\in {\cal G}$ of $G_0$.
Let \[
{\cal G}_n=\{G\in {\cal G}:D^*(Q_n,G)\in {\cal S}(Q_n)\} \]be the set of $G$'s for which $D^*(Q_n,G)$ equals a score $S_h(Q_n)$ for some uniformly bounded $h$. One can then define $G_{0n}\in {\cal G}_n$ as an approximation of $G_0$. 
\nl
{\bf $\bullet$} Let $h^*(Q_n,G_0)$ be the index so that $D^*_n(Q_n,G_{0})=S_{h^*(Q_n,G_0)}(Q_n)$.


\paragraph{Remark: Understanding ${\cal G}_n$.}
It might appear that the class of paths $\{Q_{n,\epsilon}^h:\epsilon\}$ for any bounded $h$ above is  rich enough to generate the full tangent space at $Q_n$ and thereby $D^*(Q_n,G_0)$, even for finite $n$. However, a special property of this class of paths is that it is contained in the linear span of (order $n$)   the basis functions have non-zero coefficients in $Q_n$. On the other hand, if $n$ increases, and thereby the number of basis functions converges to infinity, this class of paths will indeed be able to approximate any function in the tangent space.  Since the true $G_0$ or the relevant function of $G_0$  is generally not contained in this linear span of basis functions that make up $Q_n$,  $D^*(Q_n,G_0)\not\in {\cal S}(Q_n)$ is not contained in its set ${\cal S}(Q_n)$ of scores. For example, in the average treatment effect example, we would need that $1/\bar{G}_0(W)$ is approximated by this linear span of spline basis functions that are present in the fit $Q_n$. Therefore, indeed, there will be $G\in {\cal G}$ whose shape is such that $1/G(W)$ is in the linear span, which can then be used to define a $G_{0n}$ so that $D^*(Q_n,G_{0n})\in {\cal S}(Q_n)$. Alternatively, one directly approximates $1/\bar{G}_0(W)$ with a linear span, without being concerned if it results in a representation $1/\bar{G}_{0n}$, thereby determining an approximation $D^*_n(Q_n,G_0)$.  Since in this example $\bar{G}_0$ can be any function of $W$ with values in $(0,1)$, in this example, both methods are equivalent: i.e., if $1/\bar{G}_0$ is approximated by $\sum_j \alpha_j \phi_j$, then we can solve  for $\bar{G}_{0n}$ by setting $1/\bar{G}_{0n}=\sum_j \alpha_j\phi_j$, giving $G_{0n}=1/\sum_j \alpha_j\phi_j$.
This explains that indeed this set ${\cal G}_n$ will approximate ${\cal G}$ as  $n$ converges to infinity, so that $G_{0n}$ will approximate $G_0$, presumably certainly as fast as $Q_n$ approximates $Q_0$.
By increasing $C_n$, the number of selected basis functions in $Q_n$ will increase, thereby making the approximation $G_{0n}$ better and better.

As is evident from Theorem \ref{thfinalspline}, this approximation $G_{0n}$ should aim to approximate $G_0$ in the sense that
$R_{20}(Q_n,G_{0n},Q_0,G_0)=o_P(n^{-1/2})$ while also arranging $P_0\{D^*(Q_n,G_{0n})-D^*(Q_0,G_0)\}^2\rightarrow_p 0$.

{\bf Convenient notation for representation of $Q_n$:}
Due to finite support condition $Q^m\ll^* \mu_n$ in the definition of $m$-th order Spline HAL-MLE, we have
\begin{eqnarray}
Q_n(x)&=&\beta_n(0)+\sum_{j=0}^{m-1}\sum_{\bar{s}(j)}\beta_n(\bar{s}(j))\phi_{\bar{s}(j),\emptyset,x_s}(0)\nonumber \\
&& +\sum_{\bar{s}(m)}\sum_j \beta_n(\bar{s}(m),u_{s_m,j})\phi_{\bar{s}(m),x_s}(u_{s_m,j}).\label{splinerepresentation}
\end{eqnarray}
Recall our notation $\phi_{s,j}$ with $s\subset \{1,\ldots,k\}$ and $j\in {\cal J}_n^m(s)$, where ${\cal J}_n^m(s)$ is the finite subset of ${\cal J}^m(s)$ implied by the support points  $u_{s_m,j}$. 
Let $x_{s,j}$ be the vector of knot points (one for each component in $s$) corresponding with this basis function $\phi_{s,j}$, and we note that $\phi_{s,j}(X)=I(X(s)\geq x_{s,j})\phi_{s,j}(X)$: i.e., the support of $\phi_{s,j}$ is limited to all $x$-values for which $x(s)\geq x_{s,j}$ (and it only depends on $x$ through $x(s)$).
Analogue to (\ref{reprparspace}), we have the following representation for the $m$-th order HAL-MLE
\begin{equation}\label{convenientrepresentation}
Q_n=\sum_{s,j\in {\cal J}_n(s)}\beta_n(s,j)\phi_{s,j},\end{equation}
where we know that $\phi_{s,j}$ has support $\{x(s):x(s)\geq x_{s,j}\}$ for knot point $x_{s,j}$.

The following theorem establishes an undersmoothing condition (\ref{assumptiona}) on $C_n$ that guarantees $P_n D^*_n(Q_n,G_0)=o_P(n^{-1/2})$.
\begin{theorem}\label{thscoreequation}
Consider an approximation $D^*_n(Q_n,G_0)\in {\cal S}(Q_n)$ (i.e., scores of submodels not enforcing $L_1$-norm constant of HAL-MLE) of $D^*(Q_n,G_0)$ as defined above, and let $h_n^*$ be so that $D^*_n(Q_n,G_0)=S_{h_n^*}(Q_n)$.
Consider the representation (\ref{convenientrepresentation}) of $Q_n$. 
Note that $\beta_n$  minimizes $\beta\rightarrow P_n L\left (\sum_{s,j\in {\cal J}_n(s)}\beta_n(s,j)\phi_{s,j}\right) $ over all $\beta=(\beta(s,j):s,j\in {\cal J}_n(s))$ with $\sum_{s,j\in {\cal J}_n(s)}\mid \beta(s,j)\mid \leq C_n$. This theorem applies to any $Q_n=\sum_{s,j\in {\cal J}_n(s)}\beta_n(s,j)\phi_{s,j}$ with $\beta_n$ a minimizer of the latter empirical risk.


Assume $\pl h^*_n\pl_{\infty}=O_P(1)$,
and
\begin{equation}\label{assumptiona}
\min_{s,j\in {\cal J}_n(s),\beta_n(s,j)\not =0}\pl P_n \frac{d}{dQ_n}L(Q_n)(\phi_{s,j})\pl =o_P(n^{-1/2}).\end{equation} 
 Then, 
 \[
 P_n D^*_n(Q_n,G_0)=o_P(n^{-1.2}).\]
 
Let  $(s^*,j^*)=\arg\min_{s,j\in {\cal J}_n(s),\beta_n(s,j)\not =0}P_0\phi_{s,j}$.
 We can  replace (\ref{assumptiona}) by the following: 
 Suppose that $P_0 S_{s^*,j^*}(Q_n)^2\rightarrow_p 0$ (which will generally hold whenever $P_0\phi_{s^*,j^*}=o_P(1)$); 
 $\{S_{s,j}(Q):Q\in {\cal Q},(s,j)\}$ is contained in a Donsker class (e.g., the class of cadlag functions with uniformly bounded sectional variation norm); 
 \begin{equation}\label{suffassumptiona}
\pl P_0\left\{\frac{d}{dQ_n}L(Q_n)(\phi_{s^*,j^*})-\frac{d}{dQ_0}L(Q_0)(\phi_{s^*,j^*})\right\}\pl =o_P(n^{-1/2}),\end{equation} and $P_0 \{\frac{d}{dQ_n}L(Q_n)(\phi_{s^*,j^*}) \}^2\rightarrow_p 0$.

If we have
\[\pl P_0\left\{\frac{d}{dQ_n}L(Q_n)(\phi_{s^*,j^*})-\frac{d}{dQ_0}L(Q_0)(\phi_{s^*,j^*})\right\}
=O_P\left( P_0^{1/2}\phi_{s^*,j^*} d_0^{1/2}(Q_n,Q_0)\right);\]
 $P_0\frac{d}{dQ_n}L(Q_n)(\phi_{s^*,j^*})=O_P(P_0\phi_{s^*,j^*})$;
and  $d_0(Q_n,Q_0)=O_P(n^{-1/2-\alpha(k_1)})$ (e.g., as we showed for HAL-MLE with $\alpha(k_1)\equiv 1/\{2(k_1+2)\}$),
then (\ref{assumptiona}) can be replaced by  
\begin{equation}\label{assumptiona1}
\min_{s,j\in {\cal J}_n(s)}P_0\phi_{s,j}=o_P(n^{-1/2+\alpha(k_1)}).
\end{equation}
 \end{theorem}
 Condition (\ref{assumptiona}) is directly verifiable on the data and can thus be used to select the $m$-th order sectional variation norm bound $C_n$ for the $m$-th order Spline-HAL-MLE. For example, one could select a constant $a$ and set $C$ to the smallest value (larger than the cross-validation selector) for which the left-hand side is smaller than $a/(\sqrt{n}\log n)$ for some constant $a$. 
 The sufficient assumption (\ref{suffassumptiona}) provides understanding of what it requires in terms of $Q_n$ and $P_0$. We note that $P_0 \{\frac{d}{dQ_n}L(Q_n)(\phi_{s^*,j^*}) \}^2\rightarrow_p 0$ is a very weak condition generally implied by the support of $\phi_{s^*,j^*}$ converging to zero, and is thereby a non-condition, given our undersmoothing condition (\ref{suffassumptiona}).
The latter (\ref{suffassumptiona}) translates into the following important special case.
In this lemma we also demonstrate that if we know that $Q_n-Q_0$ converges to zero in supremum norm at a particular rate, then the support condition (\ref{assumptiona1}) can be significantly weakened. The analogue of that could have been presented in the general theorem above as well.

\begin{lemma}\label{lemmaassumptiona}
Consider the special case that $O=(Z,X)$, $L(Q)(O)$ depends on $Q$ through $Q(X)$ only, and $\frac{d}{dQ}L(Q)(\phi)=\frac{d}{dQ}L(Q)\times \phi$, i.e., the directional derivative 
$\left . \frac{d}{d\epsilon}L(Q+\epsilon \phi)\right |_{\epsilon =0}$ of $L()$ at $Q$ in the direction $\phi$ is just multiplication of a function $\frac{d}{dQ}L(Q)$ of $O$ with 
$\phi(X)$. Assume $\lim\sup_n\pl \frac{d}{dQ_n}L(Q_n)\pl_{\infty}<\infty$.
 Let $(s^*,j^*)=\arg\min_{s,j\in {\cal J}_n(s),\beta_n(s,j)\not =0}P_0\phi_{s,j}$. Assume $P_0\phi_{s^*,j^*}=o_P(1)$. Then, a sufficient condition for $P_n D^*_n(Q_n,G_0)=o_P(n^{-1/2})$ is given by (\ref{suffassumptiona}).

Assume
\[
\pl \frac{d}{dQ_n}L(Q_n)-\frac{d}{dQ_0}L(Q_0)\pl_{\infty}=O(\pl Q_n-Q_0\pl_{\infty}).\]
Then,  $P_n D^*_n(Q_n,G_0)=o_P(n^{-1/2})$ if 
\begin{equation}\label{sparsebasiscond}
\pl Q_n-Q_0\pl_{\infty} \min_{s,j\in {\cal J}_n(s),\beta_n(s,j)\not =0}P_0\phi_{s,j}=o_P(n^{-1/2}).\end{equation}
The condition (\ref{sparsebasiscond}) can be replaced by 
\[
\pl Q_n-Q_0\pl_{\infty}\min_{s,j\in {\cal J}_n(s),\beta_n(s,j)\not =0}P_n \phi_{s,j}=o_P(n^{-1/2}).\]
Here $P_0\phi_{s,j}$ and $P_n \phi_{s,j}$ can be bounded by $P_0(X(s)\geq x_{s,j}))$ and $P_n (X(s)\geq x_{s,j})$, respectively.
\end{lemma}
In \citep{vanderLaan&Bibaut17} we proved that $\pl Q_n-Q_0\pl_{\infty}\rightarrow_p 0$
under a weak absolute continuity condition, so that a  sufficient condition for (\ref{sparsebasiscond}) is $\min_{(s,j)\in {\cal J}_n,\beta_n(s,j)\not =0}P_0\phi_{s,j}=O_P(n^{-1/2})$.
However, we expect (if $m_0\geq 1$) that the rate of convergence w.r.t. supremum norm to be  $n^{-1/4}$ as achieved w.r.t. $d_0^{1/2}(Q_n,Q_0)$, in which case this only requires that $\min_{s,j\in {\cal J}_n(s),\beta_n(s,j)\not =0}P_n \phi_{s,j}=O_P(n^{-1/4})$.
The above Lemma can also be straightforwardly tailored to exploiting convergence of $d_0(Q_n,Q_0)$, as in Theorem \ref{thscoreequation}, instead of this supremum norm convergence.

\subsection{Efficiency of the Spline-HAL MLE, by including sparse basis functions}
The typical general efficiency proof used to analyze the TMLE (e.g., \citep{vanderLaan15}) can be easily generalized to the condition that $P_n D^*_n(Q_n,G_0)=o_P(n
^{-1/2})$ for some approximation $D^*_n(Q,G)$ of the actual canonical gradient $D^*(Q,G_0)$.
This results in the following theorem. 

\begin{theorem}\label{thfinalspline}
Assume condition (\ref{assumptiona}) so that $P_n D^*_n(Q_n,G_0)=o_P(n^{-1/2})$.
Assume $M_1,M_{20}<\infty$. We have $d_0(Q_n,Q_0)=O_P(n^{-1/2-\alpha(k_1)})$.

If $D^*_n(Q_n,G_0)=D^*(Q_n,G_{0n})$, then we assume
\begin{itemize}
  \item $R_2((Q_n,G_{0n}),(Q_0,G_0))=o_P(n^{-1/2})$ and 
$P_0\{D^*(Q_n,G_{0n})-D^*(Q_0,G_0)\}^2\rightarrow_p 0$.
  \item  $\{D^*(Q,G):Q\in {\cal Q},G\in {\cal G}\}$ is contained in the class of $k_1$-variate cadlag functions on a cube $[0,\tau_o ]\subset\openr^{k_1}$ in a Euclidean space and that 
$\sup_{Q\in {\cal Q},G\in {\cal G}}\pl D^*(Q,G)\pl_v^*<\infty$.  
\end{itemize}
Otherwise, we  assume
\begin{itemize}
  \item $R_2((Q_n,G_{0}),(Q_0,G_0))=o_P(n^{-1/2})$, $P_0\{D^*_n(Q_n,G_0)-D^*(Q_n,G_0)\}=o_P(n^{-1/2})$, and 
$P_0\{D^*_n(Q_n,G_{0})-D^*(Q_0,G_0)\}^2\rightarrow_p 0$.
  \item  $\{D^*_n(Q,G_0),D^*(Q,G_0):Q\in {\cal Q}\}$ is contained in the class of $k_1$-variate cadlag functions on a cube $[0,\tau_o ]\subset\openr^{k_1}$ in a Euclidean space and that 
$\sup_{Q\in {\cal Q}}\max(\pl D^*(Q,G_0)\pl_v^*,\pl D^*_n(Q,G_0)\pl_v^*)<\infty$.  
\end{itemize}
Then, $\Psi(Q_n)$ is asymptotically efficient. 
\end{theorem}
The proof is straightforward, analogue to typical efficiency proof for TMLE, and is presented in the Appendix. Regarding the condition, $P_0\{D^*_n(Q_n,G_0)-D^*(Q_n,G_0)\}=o_P(n^{-1/2})$, we note the following.
Since, for a typical choice $D^*_n(Q_n,G_0)$ in the set of scores ${\cal S}(Q_n)$, we have  $P_0 D^*_n(Q_0,G_0)=0$, so that $P_0\{D^*_n(Q_n,G_0)-D^*(Q_n,G_0)\}=P_0\{D^*_n(Q_n,G_0)-D^*_n(Q_0,G_0)\}-P_0\{D^*(Q_n,G_0)-D^*(Q_0,G_0)\}$ is indeed a second order remainder involving a product of differences $Q_n-Q_0$ and $D^*_n-D^*$.

\section{Example: HAL-MLE of treatment specific mean or average treatment effect}\label{section4}

\subsection{Formulation and relevant quantities for statistical estimation problem} 
{\bf Data and statistical model:}
Let $O=(W,A,Y)\sim P_0$,  where $Y\in \{0,1\}$ and $A\in \{0,1\}$ are binary random variables. Let $(A,W)$ have support $[0,\tau]\in \openr^{k}$, 
where $A\in [0,1]$ with only support on the edges $\{0,1\}$.
 Similarly, certain components of $W$ might be discrete so that it only has a finite set of support points in its interval.  Note $O\in [0,\tau_o]=[0,\tau]\times [0,1]$, where $[0,\tau_o]$ is a cube in Euclidean space of same dimension as $(W,A,Y)$. Let $\bar{G}(W)=E_P(A\mid W)$ and $\bar{Q}(W)=E_P(Y\mid A=1,W)$.  Assume the positivity assumption  $\bar{G}_0(W)>\delta>0$ for some $\delta>0$;
$\bar{Q}_0$, $\bar{G}_0$ are cadlag functions with $\pl \bar{Q}_0\pl_v^*\leq C^u$ and $\pl \bar{G}_0\pl_v^*\leq C^u_2$ for some finite constants $C^u,C^u_2$; $\delta<\bar{Q_0}<1-\delta$ for some $\delta>0$.  This defines the statistical model ${\cal M}$ for $P_0$.

{\bf Target parameter, canonical gradient and exact second order remainder:}
Let  $\Psi:{\cal M}\rightarrow\openr$ be defined by $\Psi(P)=E_PE_P(Y\mid W,A=1)$. 
Let $\tilde{Q}=(Q_W,\bar{Q})$, where $Q_W$ is the probability distribution of $W$. Note that $\Psi(P)=\Psi(\tilde{Q})=Q_W\bar{Q}(\cdot,1)$. We have that $\Psi$ is pathwise differentiable at $P$ with canonical gradient given by
$D^*(\tilde{Q},G)=A/\bar{G}(W)(Y-\bar{Q}(W,A))+\bar{Q}(1,W)-\Psi(\tilde{Q})$. Let $L(\bar{Q})(O)=-\left\{Y\log\bar{Q}(W,A)+(1-Y)\log(1-\bar{Q}(W,A))\right\}$ be the log-likelihood loss for $\bar{Q}$, and note that by the above bounding assumptions on $\bar{Q}$, we have that this loss function  has finite universal bounds $M_1<\infty$ and $M_{20}<\infty$.
Let $D^*_1(\bar{Q},\bar{G})=A/\bar{G}(Y-\bar{Q})$ be the $\bar{Q}$-component of the canonical gradient, $D^*_2(\tilde{Q})=\bar{Q}(1,W)-\Psi(Q)$ the $Q_W$-component, and note that
$D^*(\tilde{Q},G)=D^*_1(\bar{Q},G)+D^*_2(\tilde{Q})$. We have $\Psi(\tilde{Q})-\Psi(\tilde{Q}_0)=-P_0D^*(\tilde{Q},G)+R_{20}(\bar{Q},\bar{G},\bar{Q}_0,\bar{G}_0)$, where
\[
R_{20}(\bar{Q},\bar{G},\bar{Q}_0,\bar{Q}_0)=P_0\frac{\bar{G}-\bar{G}_0}{\bar{G} }(\bar{Q}-\bar{Q}_0).\]
{\bf Bounds on sectional variation norm and exact second order remainder:}
We have $\sup_{P\in {\cal M}}\pl D^*(\tilde{Q}(P),G(P))\pl_v^*<C(C^u,C^u_2)$ for some finite constant $C$ implied by the universal bounds $(C^u,C^u_2)$ on the sectional variation norm of $\bar{Q},\bar{G}$. We also note that, using Cauchy-Schwarz inequality,  $R_{20}(\bar{Q},\bar{G},\bar{Q}_0,\bar{G}_0)\leq \frac{1}{\delta}\pl \bar{Q}-\bar{Q}_0\pl_{P_0}\pl \bar{G}-\bar{G}_0\pl_{P_0}$, where
$\pl f\pl_{P_0}^2=\int f^2(o)dP_0(o)$. 

\subsection{HAL-MLE}
Let $Q=\mbox{Logit}\bar{Q}$ and let's $L(Q)(O)=-A\{Y\log \bar{Q}(W)+(1-Y)\log(1-\bar{Q}(W))\}$ be the log-likelihood loss restricted to the observations with $A=1$. Let ${Q}_{C,n}=\arg\min_{{Q},\pl {Q}\pl_v^*<C}P_n L(Q)$ be the $C$-specific $0$-order Spline-HAL-MLE for a given bound $C$ on the sectional variation norm. 
Let $C_n\leq C^u$ be a data adaptive selector that is larger or equal than the cross-validation selector, so that $P(C_{n,cv}\leq C_n\leq C^u)=1$. Let ${Q}_n={Q}_{C_n,n}$, and $Q_{W,n}$ be the empirical probability measure of $W_1,\ldots,W_n$. 
We can represent $Q_n=\sum_{(s,j\in \in {\cal J}_n(s)}\beta_n(s,j)\phi_{s,j}$, where $\phi_{s,j}=I(W(s)\geq w_{s,j})$ for a knot point $w_{s,j}$ varying over all observations $\{W_{s,i}:i=1,\ldots,n\}$ across all subsets $s\subset\{1,\ldots,k_1\}$.
By our rate of convergence results on the HAL-MLE we have that $\pl {Q}_n-{Q}_0\pl_{P_0}=O_P(n^{-1/4-\alpha(k_1)})$ for $\alpha(k)=1/2(k+2)$.
The HAL-MLE of $\Psi(\tilde{Q}_0)$ is the plug-in estimator $\Psi(\tilde{Q}_n)=Q_{W,n}\bar{Q}_n$. Note that $P_n D^*_2(\tilde{Q}_n)=0$ for any ${Q}_n$.
Thus, we are only concerned with showing that $P_n D^*_1(Q_n,G_{0})=o_P(n^{-1/2})$.

{\bf Class of paths absolute continuous w.r.t. $Q_n$:}
Consider the following class of paths 
 \[
{Q}_{n,\epsilon}^h(x)={Q}_n+\epsilon f(h,{Q}_n),\]
where
\[
f(h,{Q}_n)(x)=h(0){Q}_n(0)+\sum_{s\subset \{1,\ldots,m\}}\int_{(0_s,x_s]}  h(s,u_s)) {Q}_{n,s}(du_s).\]
This defines a path $\{{Q}_{n,\epsilon}^h:\epsilon\in (-\delta,\delta)\}$ for each uniformly bounded function $h$, as in our general representation.
 
 {\bf Set of scores generated by class of paths:}
  The scores generated by this family of paths are given by:
 \[
S_h({Q}_n)\equiv \left .  \frac{d}{d\epsilon}L({Q}_{n,\epsilon}^h)\right |_{\epsilon =0}  =A f(h,{Q}_n)(Y-{Q}_n(W)).\]
This defines a set of scores ${\cal S}(Q_n)=\{S_h(\bar{Q}_n): \pl h\pl_{\infty}<\infty\}$. 
Note that in order to solve for an $h$ so that $S_h(\bar{Q}_n)=D^*_1(\bar{Q}_n,\bar{G}_0)$ would require $f(h,\bar{Q}_n)(W)=1/\bar{G}_0(W)$. However, since $\bar{G}_0$ is not sectional absolute continuous w.r.t. ${Q}_n$ ( i.e., ${Q}_{n,s}$ is discrete for all subsets $s$, while $\bar{G}_{0,s}$ is (say) continuous), there does not exist a $h$ for which $f(h,{Q}_n)=1/\bar{G}_0$.
Thus, $D^*({Q}_n,\bar{G}_0)\not\in \{S_h({Q}_n):\pl h\pl_{\infty}<\infty\}$.

{\bf Score equations solved by HAL-MLE:}\newline
Let $r(h,{Q}_n)\equiv h(0))\mid {Q}_n(0)\mid +\sum_{s\subset \{1,\ldots,m\}}\int_{(0_s,x_s]} h(s,u_s))\mid d{Q}_{n,s}(u_s)\mid$.
The HAL-MLE solves
\[
P_n S_h({Q}_n)=0\mbox{ for all $h$ with $r(h,{Q}_n)=0$.}\]

\subsection{Defining approximation $G_{0n}$}
We define
\[
{\cal G}_n\equiv \{\bar{G}\in {\cal G}: \bar{G}\ll^* \bar{Q}_n\}.\]
We note that if $\bar{G}_{s}\ll \bar{Q}_{n,s}$, then we also have $1/\bar{G}_s\ll \bar{Q}_{n,s}$ as well. Here we use that if $g(x)=1/f(x)$, then $g_s(dx_s)=-1/f^2_s(x_s) f_s(dx_s)$.
Therefore, if $\bar{G}\ll^*\bar{Q}_n$, then  we can find a $h$ so that 
$f(h,{Q}_n)(A,W)=A/\bar{G}(W)$, and thereby that $D^*_1({Q}_n,\bar{G})=S_h(Q_n)$.

Let 
\[
G_{0n}=\arg\min_{\bar{G}\in {\cal G}_n}\pl \bar{G}-\bar{G}_0\pl_{P_0},\]
where $\pl \bar{G}-\bar{G}_0\pl_{P_0}$ is the $L^2(P_0)$-norm of $\bar{G}-\bar{G}_0$. 
Then, $D^*_1({Q}_n,\bar{G}_{0n})\in \{S_h({Q}_n):h\}$ so that we can find a $h^*({Q}_n,\bar{G}_0)$ so that
\[
D^*_1({Q}_n,\bar{G}_{0n})=S_{h^*({Q}_n,\bar{G}_0)}({Q}_n).\]

\subsection{Application of Theorem \ref{thfinalspline}}
We need to assume  $R_2((Q_n,G_{0n}),(Q_0,G_0))=o_P(n^{-1/2})$ and 
$P_0\{D^*(Q_n,G_{0n})-D^*(Q_0,G_0)\}^2\rightarrow_p 0$.
 The latter already holds if  $\pl \bar{G}_{0n}-\bar{G}_0\pl_{P_0}\rightarrow_p 0$. However, the first condition relies on a rate of convergence. For example, this will hold if $\pl \bar{G}_{0n}-\bar{G}_0\pl_{P_0}=O_P(n^{-1/4})$. This appears to be  a reasonable condition, since $\bar{G}_{0n}$ is the $L^2(P_0)$-projection of $\bar{G}_0$ onto ${\cal G}_n$, so that the only concern would be that the set ${\cal G}_n$ does not approximate fast enough ${\cal G}$ as $n$ converges to infinity. 
However,  if the set of basis functions  is rich enough for $\bar{Q}_n$ to converge at a rate faster than $n^{-1/4}$ to $\bar{Q}_0$ (not allowing to choose the coefficients based on $P_0$), then the resulting linear combination  of indicator basis functions should generally also be rich enough for approximating the true $G_0$ with a rate $n^{-1/4}$ (now allowing to select the coefficients of the basis functions in terms of $G_0$). 

{\bf Verification of Assumption \ref{assumptiona} of Theorem \ref{thfinalspline}:}
Assumption (\ref{assumptiona}) is stating that
\[
\begin{array}{l}
\min_{(s,j)\in {\cal J}_n} P_n \frac{d}{dQ_n}L(Q_n)(\phi_{s,j})\\
 =
2\frac{1}{n}\sum_i \phi_{s,j}(1,W_i)I(A_i=1)I(W_{s,i}>w_{s,j})(Y_i-\bar{Q}_n(1,W_i))\\
=o_P(n^{-1/2}).
\end{array}
\]
We apply the last part of Theorem \ref{thscoreequation}.
Since $\frac{d}{dQ}L(Q)(\phi)=\phi(A,W)(Y-\bar{Q}(A,W))$, it follows that have
\begin{equation}\label{helpz1}
\pl \frac{d}{dQ_n}L(Q_n)-\frac{d}{dQ_0}L(Q_0)\pl_{P_0}=O(\pl Q_n-Q_0\pl_{P_0}).\end{equation}
Given that we have $d_0(Q_n,Q_0)=O_P(n^{-1/2-\alpha(k_1)})$, it follows that the remaining condition is (\ref{assumptiona1}), or, equivalently,
\[
\min_{s,j\in {\cal J}_n(s),\beta_n(s,j)\not =0}P_n\phi_{s,j}=O_P(n^{-1/2+\alpha(k_1)}).\]
This reduces to the assumption that 
 $O\left( \min_{\{s,j\in {\cal J}_n(s):\beta_n(s,j)\not =0\}} P_n(W(s)\geq  w_{s,j})\right)=O_P(n^{-1/2+\alpha(k_1)})$. We arrange this assumption to hold by selecting $C_n$ accordingly. 
 Similarly, we can apply Lemma \ref{lemmaassumptiona} but now expressing the latter condition in terms of $\pl Q_n-Q_0\pl_{\infty}$.

This proves the following efficiency theorem  for the HAL-MLE in this particular estimation problem.
\begin{theorem}\label{thexample}
Consider the formulation above of the statistical estimation problem. 
Let
\[
{\cal G}_n=\{\bar{G}\in {\cal G}: \bar{G}\ll^* \bar{Q}_n\},\] and
\[
\bar{G}_{0n}=\arg\min_{\bar{G}\in {\cal G}_{n}} \pl G-G_0\pl_{P_0}.\]

{\bf Assumptions:}
\begin{itemize}
\item 
$\pl \bar{G}_{0n}-\bar{G}_0\pl_{P_0}=O_P(n^{-1/4})$, where we can use that $\pl Q_n-Q_0\pl_{P_0}=o_P(n^{-1/4})$.
\item 
Given the  fit ${Q}_n=\sum_{s,j\in {\cal J}_n(s)}\beta_n(s,j)\phi_{s,j}$ with support points the observations $\{W_j(s): j=1,\ldots,n,s\}$ and indicator basis functions $\phi_{s,j}(W)=I(W(s)>W_j(s))$, we assume that $C_n<C^u$ for some finite $C^u$ is chosen so that either
\[
\pl Q_n-Q_0\pl_{\infty}\min_{\{s,j\in {\cal J}_n(s):\beta_n(s,j)\not =0\}}  P_n(W(s)\geq W_j(s))=o_P(n^{-1/2}),\]
or 
\[
\min_{s,j\in {\cal J}_n(s),\beta_n(s,j)\not =0}P_n\phi_{s,j}=O_P(n^{-1/2+\alpha(k_1)}).\]
\end{itemize}
Then, $\Psi(Q_n)$ is an asymptotically efficient estimator of $\Psi(Q_0)$. 

\end{theorem}

\section{Example: HAL-MLE for the integral of the square of the data density}
\label{section5}
Let $O\sim P_0$ be a $k_1$-variate random variable with Lebesgue density $p_0$ that is assumed to be bounded from below by a $\delta>0$ and from above by an $M<\infty$.
Let $\{P_Q:Q\in {\cal Q}\}$ be a parametrization of the probability measure of $O$   in terms of a functional parameter $Q$ that varies over a class of  multivariate real valued cadlag functions on $[0,\tau]$ with finite sectional variation norm. Below we will focus on the particular parameterization given by  $p_Q=c(Q)\{\delta+(M-\delta)\mbox{expit}(Q)\}$, where
$\mbox{expit}(x)=1/(1+\exp(-x))$, and $c(Q)$ is the normalizing constant defined by $\int p_Qdo=1$.
Note that in this parameterization $Q$ can be any cadlag function with finite sectional variation norm, thereby allowing that the densities $p_Q$ are discontinuous (but cadlag). 
Another possible parametrization is obtained through the following steps:  1) modeling the density $p(x)$ as a product $\prod_{j=1}^k p_j(x_j\mid \bar{x}(j-1))$ of univariate conditional densities of $x_j$, given $\bar{x}(j-1)$; 2) modeling each univariate conditional density $p_j$ in terms of its univariate conditional hazard $\lambda_j$; 3) modeling this hazard as $\lambda(x_j\mid \bar{x}(j-1))= \exp(Q_j(x_j,\bar{x}(j-1)))$ (or discretizing it and modeling it with a logistic function in $Q_j$), and 4) setting 
$Q=(Q_1,\ldots,Q_{k_1})$.
With this latter parametrization each $Q_j$ varies over a parameter space of cadlag functions with finite sectional variation norm. 

Let the statistical model ${\cal M}=\{P_Q: Q\in {\cal Q}(C^u)\}$ for $P_0$ be nonparametric beyond that each probability distribution is dominated by the Lebesgue measure, $Q$ varies over cadlag functions with sectional variation norm bounded by $C^u$. 
The statistical target parameter $\Psi:{\cal M}\rightarrow\openr$ is defined by $\Psi(P)=\int p^2(o)do$. The canonical gradient of $\Psi$ at $P$ is given by $D^*(P)(O)=2 (p(O)-\Psi(P))$, and, the exact second order remainder $R_2(P,P_0)=\Psi(P)-\Psi(P_0)+P_0D^*(P)$ is given by $R_2(P,P_0)=-\int (p-p_0)^2(o) do$. 

Let $L(Q)=-\log p_Q$ be the log-likelihood loss function for $Q$. 
Let $Q_n$ be a $m=0$-order HAL-MLE bounding the sectional variation norm by a $C_n<C^u$.
We wish to establish conditions on $C_n$ so that $\Psi(Q_n)=\int p_{Q_n}^2do$ is an asymptotically efficient estimator of $\Psi(Q_0)=\int p_{Q_0}^2 d0$.
We assume this HAL-MLE is discrete so that we can use the finite dimensional representation $Q_n=\sum_{s,j\in {\cal J}_n(s)}\beta_n(s,j)\phi_{s,j}$ with $\pl \beta_n\pl_{L_1}\leq C_n$, as in our general presentation. Let $Q_{n,\epsilon}^h=Q_n(0)(1+\epsilon h(0))+\sum_{s}\int_{(0_s,x_s]} (1+\epsilon h(s,u_s))dQ_{n,s}(u_s)$, indexed by any bounded function $h$, be the paths as defined in our general presentation (and previous section). Let $S_h(Q_n)=\left . \frac{d}{d\epsilon}L(Q_{n,\epsilon}^h)\right |_{\epsilon =0}$ be score of this path under the log-likelihood loss. 
These scores are given by \[
S_h(Q_n)=\frac{d}{dQ_n}L(Q_n)(f(h,Q_n)),\]
 where $f(h,Q_n)=Q_n(0)h(0)+\sum_{s}\int_{(0_s,x_s]} h(s,u_s)dQ_{n,s}(u_s)$.  Let ${\cal S}(Q_n)=\{S_h(Q_n):h\}$ be the collection of scores.
 In order to apply Theorem \ref{thfinalspline} we need to determine an approximation $D_n(Q_n)\in {\cal S}(Q_n)$ of the canonical gradient $D^*(Q_n)=2(p_{Q_n}-\Psi(Q_n))$.
We have
\[
\begin{array}{l}
S_h(Q)=A(f(h,Q))/C(Q)-M\frac{\exp(Q)}{(1+\exp(Q))(\delta+\delta\exp(Q)+M)} f(h,Q),
\end{array}
\]
where
\[A(f)=\frac{\int \exp(Q)/(1+\exp(Q))^2 f d0}{(\int (\delta+M/(1+\exp(Q))) d0)^2}.\]
Let $G(Q)=-M\frac{\exp(Q)}{(1+\exp(Q))(\delta+\delta\exp(Q)+M)}$, so that the equation
$S_h(Q)=D^*(Q)$ corresponds with $G(Q)f(h,Q)+C(Q)^{-1}A(f(h,Q))=D^*(Q)$, which can be rewritten as
$f(h,Q)+G_1(Q)A(f(h,Q))=D^*(Q)/G(Q)$, and $G_1(Q)=1/(C(Q)G(Q))$.
Let $D_1(Q)=D^*(Q)/G(Q)$, so that the equation becomes $f+G_1(Q) A(f)=D_1(Q)$. Once we have solved for $f$, whose solution we will denote with $f(Q)$, then it remains to solve
for $h$ in $f(h,Q)=f(Q)$ or find a closest solution. 
It is important to note the $f\rightarrow A(f)$ is a linear real valued operator. Applying this operator to both sides yields
$A(f)+A(f) A(G_1(Q))=A(D_1(Q))$, so that we obtain the solution
\[
A(f)=\frac{ A(D_1(Q))}{1+A(G_1(Q))}.\]
Plugging this back into the equation, we obtain $f(Q)\equiv D_1(Q)-\frac{ G_1(Q)A(D_1(Q))}{1+A(G_1(Q))}$. 
Thus, we have shown that if we can set $f(h,Q_n)=f(Q_n)$, then we have $S_h(Q_n)=D^*(Q_n)$.
It remains to determine  a choice $h(Q_n)$ so that $f(h,Q_n)\approx f(Q_n)$.
The space $\{f(h,Q_n):h\}$ equals $\{\sum_{s,j\in {\cal J}_n(s)}\alpha(s,j)\phi_{s,j}:\alpha\}$ the linear span of the basis functions $\{\phi_{s,j}:s,j\in {\cal J}_n(s)\}$.
Let $f_n(Q_n)$ be the projection of $f(Q_n)$ onto this linear space, for example, w.r.t. $L^2(P_0)$-norm.
Let $h_n(Q_n)$ be the solution of $f(h,Q_n)=f_n(Q_n)$, and let $D^*_n(Q_n)=S_{h_n(Q_n)}(Q_n)$ be our desired approximation of $D^*(Q_n)$ which is an element of the set of scores $\{S_h(Q_n):h\}$.
We note that 
\[
\begin{array}{l}
D^*_n(Q_n)-D^*(Q_n)=S_{h_n(Q_n)}(Q_n)-D^*(Q_n)\\
=
G(Q_n)f_n(Q_n)+C(Q_n)^{-1}A(f_n(Q_n))-D^*(Q_n)\\
=
G(Q_n)f_n(Q_n)+C(Q_n)^{-1}A(f_n(Q_n))-G(Q_n)f(Q_n)-C(Q_n)^{-1}A(f(Q_n))\\
=G(Q_n)(f_n(Q_n)-f(Q_n))+C(Q_n)^{-1}A(f_n(Q_n)-f(Q_n)).
\end{array}
\]
We will assume that $\pl f_n(Q_n)-f(Q_n)\pl_{P_0}=o_P(n^{-1/4})$.
The main condition beyond (\ref{assumptiona}) of Theorem \ref{thfinalspline} is that
$P_0\{D^*_n(Q_n)-D^*(Q_n)\}=o_P(n^{-1/2})$.
Note that $P_0 D^*_n(Q_0)=0=P_0D^*(Q_0)$. Therefore, 
\[
\begin{array}{l}
P_0\{D^*_n(Q_n)-D^*(Q_n)\}=P_0\{D^*_n(Q_n)-D^*_n(Q_0)\}-P_0\{D^*(Q_n)-D^*(Q_0)\}\\
=P_0\{
G(Q_n)(f_n(Q_n)-f(Q_n))\}+P_0\{C(Q_n)^{-1}A(f_n(Q_n)-f(Q_n))\}\\
-P_0\{G(Q_0)(f_n(Q_0)-f(Q_0))-C(Q_0)^{-1}A(f_n(Q_0)-f(Q_0)) \}
.
\end{array}
\]
Let $\Pi_n$ be the projection operator on the linear span generated by the basis function of $Q_n$, which is of the same dimension as the number of basis functions.
The latter difference can also be represented as
 \[
 P_0\{ D^*(Q_n)-D^*(Q_0)-\Pi_n(D^*(Q_n)-D^*(Q_0))\},\]
 or, if we define $\Pi_n^{\perp}=(I-\Pi_n)$ as the projection operator onto the orthgonal complement of the linear space spanned by the basis functions in $Q_n$, then this term can be denoted as
 \begin{equation}\label{keycondex}
 P_0\{ \Pi_n^{\perp}(D^*(Q_n)-D^*(Q_0))\},\end{equation}
 which can, in particular, be bounded by the operator norm $\pl \Pi_n^{\perp}\pl$ of $\Pi_n^{\perp} $ times the $L^2(P_0)$-norm of $D^*(Q_n)-D^*(Q_0)$. Thus, if we assume that $\pl \Pi_n^{\perp}\pl=O_P(n^{-1/4})$, then it follows that this term is $o_P(n^{-1/2})$.
 We will simply assume (\ref{keycondex}) to be $o_P(n^{-1/2})$.
The other conditions, beyond (\ref{assumptiona}) of Theorem \ref{thfinalspline} hold by the fact that $\pl Q_n-Q_0\pl_{P_0}=o_P(n^{-1/4})$ and that $D^*(Q_n),D^*_n(Q_n)$ falls in a $P_0$-Donsker class of cadlag functions with universal bound on sectional variation norm. 

{\bf Verification of Assumption \ref{assumptiona} of Theorem \ref{thfinalspline}:}
Assumption (\ref{assumptiona}) is stating that
\[
\begin{array}{l}
\min_{s,j\in {\cal J}_n(s)} P_n \frac{d}{dQ_n}L(Q_n)(\phi_{s,j})
=o_P(n^{-1/2}).
\end{array}
\]
We apply the last part of Theorem \ref{thscoreequation}.
We have 
\begin{equation}\label{helpz2}
\pl \frac{d}{dQ_n}L(Q_n)-\frac{d}{dQ_0}L(Q_0)\pl_{P_0}=O(\pl Q_n-Q_0\pl_{P_0}).\end{equation}
Given that we have $d_0(Q_n,Q_0)=O_P(n^{-1/2-\alpha(k_1)})$, it follows that the remaining condition is (\ref{assumptiona1}), or, equivalently,
\[
\min_{s,j\in {\cal J}_n(s),\beta_n(s,j)\not =0}P_n\phi_{s,j}=O_P(n^{-1/2+\alpha(k_1)}).\]
This reduces to the assumption that 
 $O\left( \min_{\{s,j\in {\cal J}_n(s):\beta_n(s,j)\not =0\}} P_n(O(s)\geq  O_{s,j})\right)=O_P(n^{-1/2+\alpha(k_1)})$. We arrange this assumption to hold by selecting $C_n$ accordingly. 

This proves the following efficiency theorem  for the HAL-MLE in this particular estimation problem.
\begin{theorem}\label{thexample2}
Let $O\sim P_0$ be a $k_1$-variate random variable with Lebesgue density $p_0$ that is assumed to be bounded from below by a $\delta>0$ and from above by an $M<\infty$.
Let  $p_Q=c(Q)\{\delta+(M-\delta)\mbox{expit}(Q)\}$, where
$\mbox{expit}(x)=1/(1+\exp(-x))$, and $c(Q)$ is the normalizing constant defined by $\int p_Qdo=1$, where $Q\in {\cal Q}(C^u)$ can be any cadlag function with finite sectional variation norm bounded by $C^u$.
Let the statistical model ${\cal M}=\{P_Q: Q\in {\cal Q}(C^u)\}$ for $P_0$ be nonparametric beyond that each probability distribution is dominated by the Lebesgue measure, $Q$ varies over cadlag functions with sectional variation norm bounded by $C^u$. 
The statistical target parameter $\Psi:{\cal M}\rightarrow\openr$ is defined by $\Psi(P)=\int p^2(o)do$, which we also denote with $\Psi(Q)$. The canonical gradient of $\Psi$ at $P$ is given by $D^*(P)(O)=2 (p(O)-\Psi(P))$, and, the exact second order remainder $R_2(P,P_0)=\Psi(P)-\Psi(P_0)+P_0D^*(P)$ is given by $R_2(P,P_0)=-\int (p-p_0)^2 do$.

Consider the formulation above of the statistical estimation problem. 

{\bf Assumptions:}
\begin{itemize}
\item 
 $\pl Q_n-Q_0\pl_{P_0}=o_P(n^{-1/4})$.

\item
Given the  fit ${Q}_n=\sum_{s,j\in {\cal J}_n(s)}\beta_n(s,j)\phi_{s,j}$ with support points the observations $\{O_j(s): j=1,\ldots,n,s\}$ and indicator basis functions $\phi_{s,j}(W)=I(O(s)>O_j(s))$, we assume that $C_n<C^u$ for some finite $C^u$ is chosen so that either
\[
\pl Q_n-Q_0\pl_{\infty}\min_{\{s,j\in {\cal J}_n(s):\beta_n(s,j)\not =0\}}  P_n(W(s)\geq W_j(s))=o_P(n^{-1/2}),\]
or 
\[
\min_{s,j\in {\cal J}_n(s),\beta_n(s,j)\not =0}P_n\phi_{s,j}=O_P(n^{-1/2+\alpha(k_1)}).\]
\item  Let $\Pi_n^{\perp}$ be the projection operator in $L^2(P_0)$ onto the orthogonal complement of the linear span of the basis functions $\{\phi_{s,j}:s,j\in {\cal J}_n(s),\beta_n(s,j)\not =0\}$ in the fit of $Q_n$. Assume 
\begin{equation}\label{keycondexth}
 P_0\{ \Pi_n^{\perp}(D^*(Q_n)-D^*(Q_0))\}=o_P(n^{-1/2}).\end{equation}
A sufficient condition is that the operator norm $\pl \Pi_n^{\perp}\pl$ of $\Pi_n^{\perp}$ is $O_P(n^{-1/4})$.
\end{itemize}
Then, $\Psi(Q_n)$ is an asymptotically efficient estimator of $\Psi(Q_0)$. 

\end{theorem}

\section{Simulation study}\label{section6}

Our global undersmoothing condition only specifies a sufficient rate at which the sparsest selected basis function should converge to zero, but it does not provide a constant in front of this rate. Thus, it does not immediately yield a practical method for tuning the level of undersmoothing. In our simulation studies, we investigate the targeted $L_1$-norm selector that is chosen so that the empirical mean of the canonical gradient at the HAL-MLE (indexed by $L_1$-norm) and possibly a HAL-MLE of the nuisance parameter in the canonical gradient is $o_P(n^{-1/2})$. In extensive simulations, this method appears to give better practical results than several direct implementations of our global undersmoothing criterion (i.e., the choice of constant matters for practical performance). More research will be needed to investigate if one can construct a global undersmoothing selector (according to our theorem) that would result in well behaved efficient plug-in estimators across a large class of target parameters. Our simulations also demonstrate that our targeted selection method for  undersmoothing controls the sectional variation norm of the fit, which is a crucial part of the Donsker class or asymptotic equicontinuity condition. 

\subsection{Simulations for the ATE}

We simulated a vector $W = (W_1, W_2)$. $W_1$ was simulated by drawing $Z \sim \mbox{Beta}(0.85, 0.85)$ and setting $W_1 = 4Z -2$. $W_2$ was drawn independently from a Bernoulli(0.5) distribution. Given $W = w$, a binary random variable $A$ was drawn with probability $A=1$ equal to $\bar{G}_0(w) = \mbox{logit}^{-1}\{w_1 - 2w_1w_w\}$. Given $W = w$, we set $Y = \bar{Q}_0(w) + \epsilon$, where $\bar{Q}_0(w) = \mbox{logit}^{-1}\{w_1 - 2w_1w_2\}$ and $\epsilon \sim \mbox{Normal}(0, 0.25)$. We refer readers back to Section \ref{section4} for the definition of the target parameter and form of the canonical gradient. 

We built our undersmoothed estimator of $\Psi(P_0)$ as follows. We estimate $\bar{Q}_0$ using a HAL regression estimator and select the regularization of the estimator by choosing the smallest value of $C$ for $L_1$-norm such that $$
	P_n D^*(Q_{C,n}, \bar{G}_n) < \frac{P_n \{D^*(Q_{C,n}, \bar{G}_n)^2\}}{\mbox{log}(n) n^{1/2}} \ ,
$$
where $\bar{G}_n$ is the HAL-MLE estimate of $\bar{G}_0$ (i.e., a HAL regerssion that uses cross-validated choice for $C$). We then computed the plug-in estimator as described in Section \ref{section4}. 

We generated 3,000 data sets in this way and computed the undersmoothed HAL estimate. We report the estimator's bias (scaled by $n^{1/2}$), variance (scaled by $n$), mean squared error (by $n$), and the sampling distribution of $n^{1/2}\{\Psi(\tilde{Q}_n) - \Psi(P_0)\}$. We additionally report on the behavior of $n^{1/2} P_n D^*(Q_{C_n,n},\bar{G}_0)$ and $$
n^{1/2} \left\{\min_{s,j\in {\cal J}_n(s),\beta_n(s,j)\not =0}\pl P_n \frac{d}{dQ_n}L(Q_n)(\phi_{s,j})\pl \right\} \ . 
$$

\begin{figure}
\centering
\includegraphics[width=0.45\textwidth]{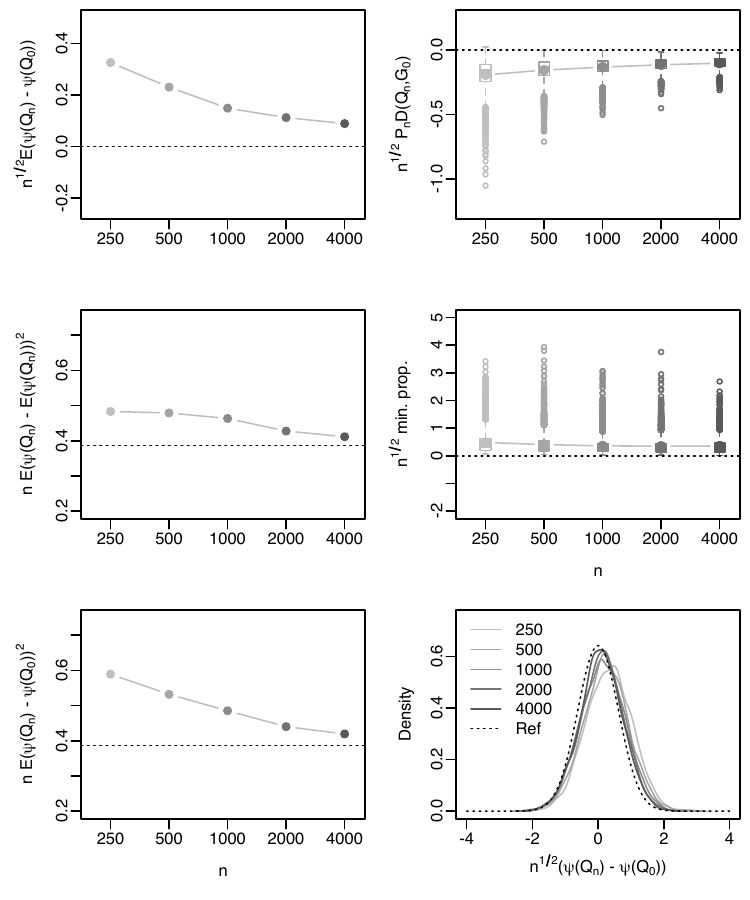}
\caption{Left column top to bottom: bias, variance, and mean squared-error (all scaled by $n^{1/2}$) of undersmoothed HAL-MLE. Right column top to bottom: scaled empirical average of canonical gradient, empirical average of quantity given in equation (\ref{assumptiona}), sampling distribution of scaled and centered estimator. The dashed lines in the variance and mean-squared error plots denote the efficiency bound. The reference sampling distribution for the estimators is a mean-zero Normal distribution with this variance.}
\label{sim_results}
\end{figure}

As predicted by theory, the bias of the estimator diminishes faster than $n^{-1/2}$ and the variance of the estimator approaches the efficiency bound in larger samples (Figure \ref{sim_results}). The empirical average of the canonical gradient is appropriately controlled (top right) and our selection criteria for the HAL tuning parameter appears to also satisfy the global criteria stipulated by equation (\ref{assumptiona}). At all sample sizes, the sampling distribution of the scaled and centered estimator is well-approximated by the efficient asymptotic distribution. 

\subsection{Simulations for the integral of the square of the density}

We simulated a univariate variable $O \sim N(-4, 5/3)$ and evaluated the performance of undersmoothed HAL for estimating the integral of the square of the density of $O$ (Section \ref{section5}). We implemented a HAL-based estimator of the density using an approach similar to the one described in \citet{munoz2011super}. This approach entails estimating a discrete hazard function using HAL using a pre-specified binning of the real line. For this simulation, we used 320 equidistant bins, and note that the HAL density estimator is robust to this choice, so long as a large enough value is chosen. We sample 1,000 data sets for each of several sample sizes ranging from $n=100$ to 100,000. We compare the results for undersmoothed HAL to those obtained by using a typical implementation of HAL that selects the level of smoothing based on cross-validation. We compared these estimators on the same criterion described in the previous subsection. 

The simulations results reflect what is expected based on theory. In particular, the undersmoothed HAL achieves the efficiency bound in large samples and the scaled-centered sampling distribution of the estimator is well approximated by the efficient asymptotic distribution. We found that our selection criterion for the level of undersmoothing based on the EIF led to control of the variation norm of the resultant fit. On the other hand, results for the HAL estimator with level of smoothing selected via cross-validation demonstrated that this estimator does not have bias that is decreasing faster than $n^{-1/2}$. Thus, this estimator performs worse in terms of all criteria that we considered. 

\begin{figure}[ht]
    \centering
    (a)\includegraphics[width=0.4\textwidth]{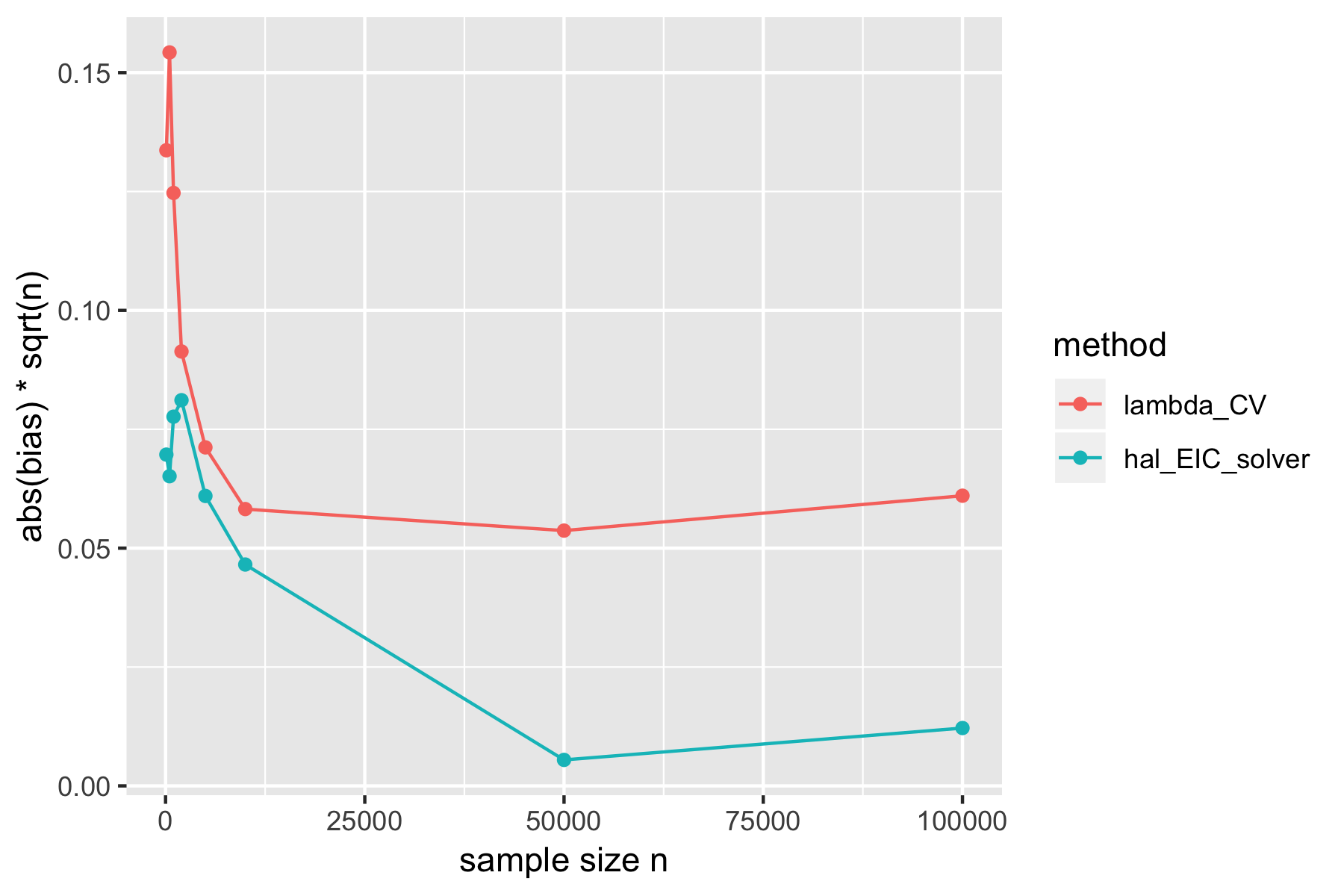}
    (b)\includegraphics[width=0.4\textwidth]{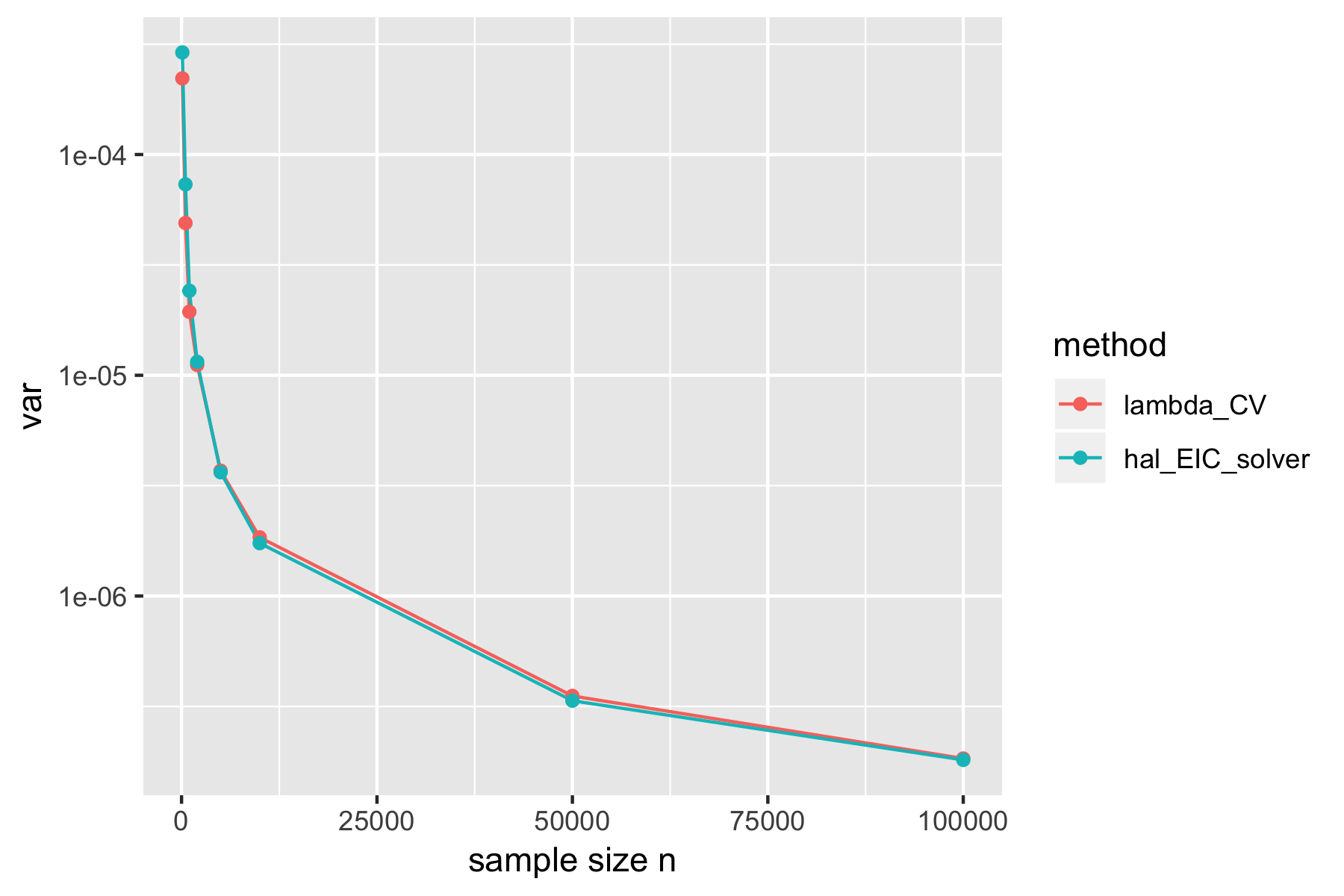}
    (c)\includegraphics[width=0.4\textwidth]{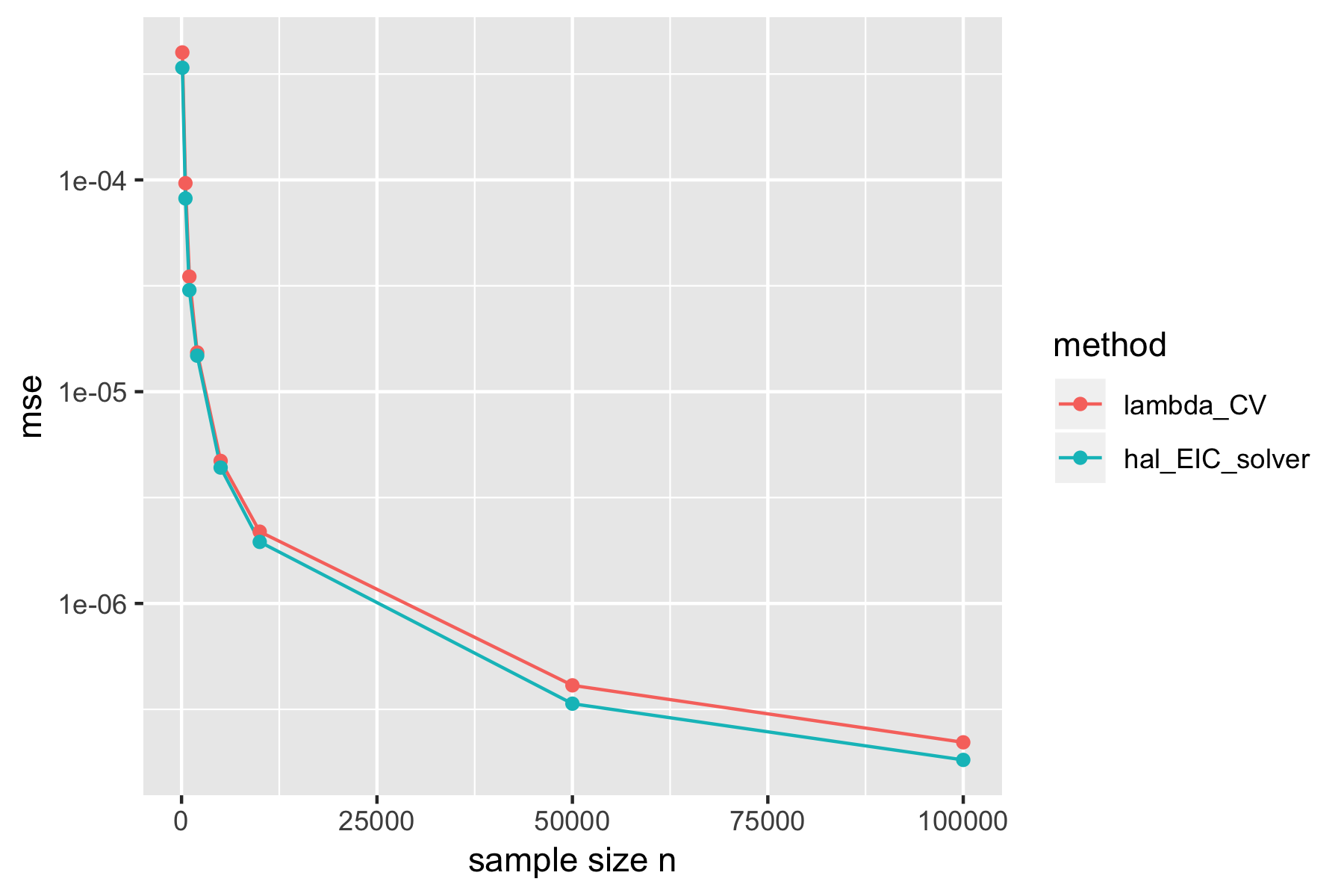}
    (d)\includegraphics[width=0.4\textwidth]{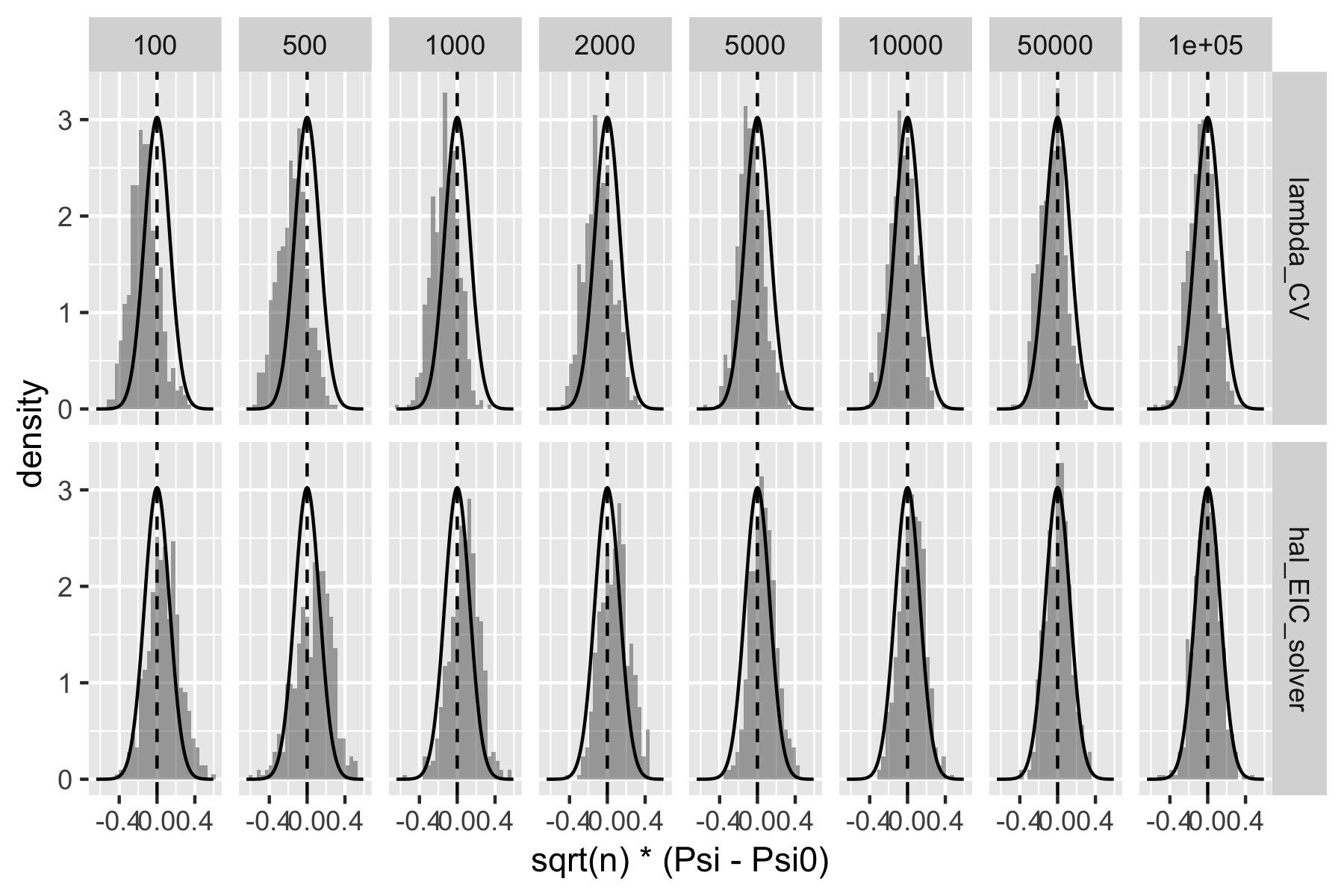}
    \caption{Simulation results for the average density value parameter: (a) $\sqrt{n}$ times absolute bias (b) variance (c) MSE (d) histogram of $\sqrt{n}(\Psi_n - \Psi_0)$ for undersmoothed HAL and HAL(CV) at different sample sizes}
    \label{fig:avgdens1}
\end{figure}

\begin{figure}[ht]
    \centering
    (a)\includegraphics[width=0.4\textwidth]{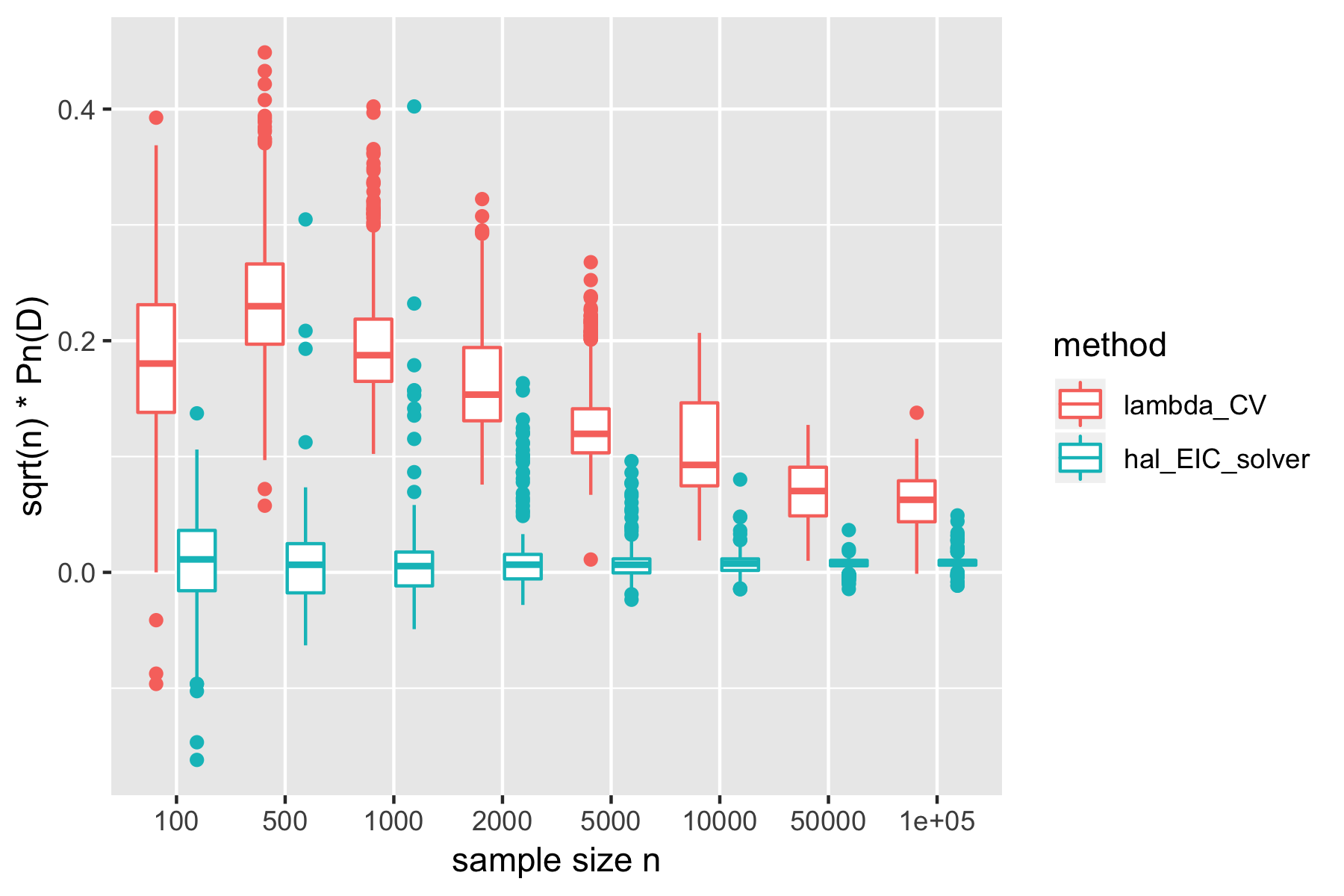}
    (b)\includegraphics[width=0.4\textwidth]{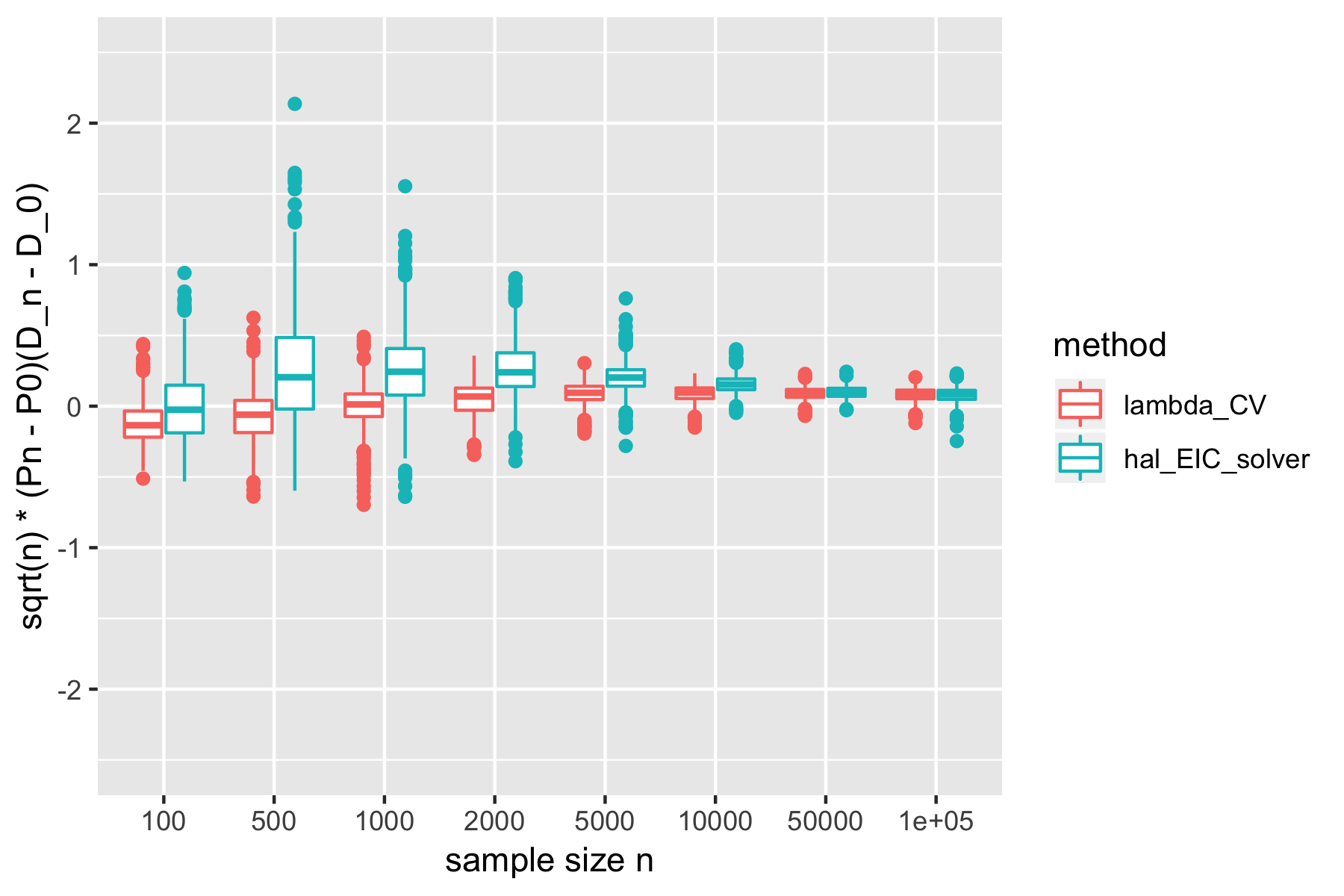}
    (c)\includegraphics[width=0.4\textwidth]{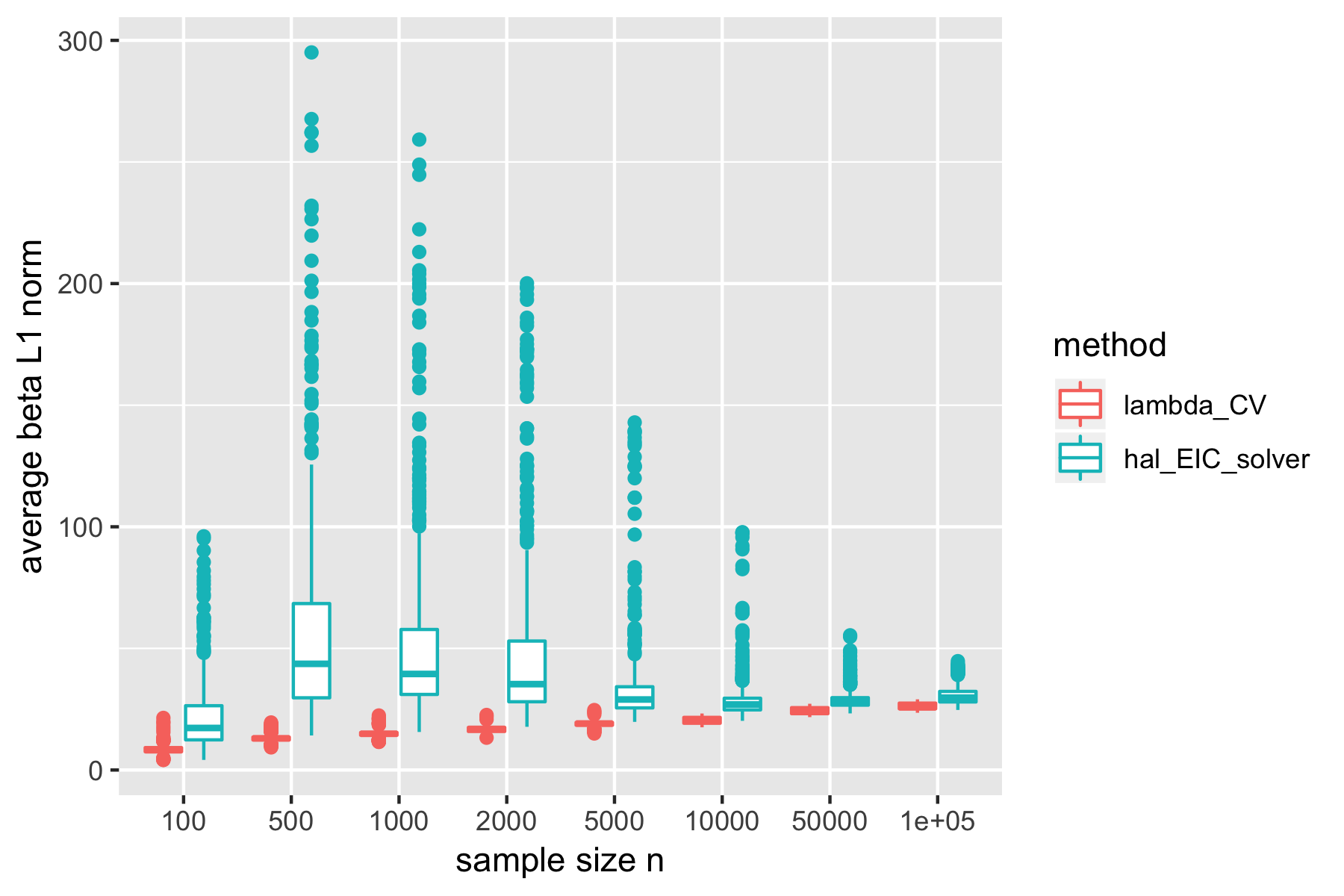}
    (d)\includegraphics[width=0.4\textwidth]{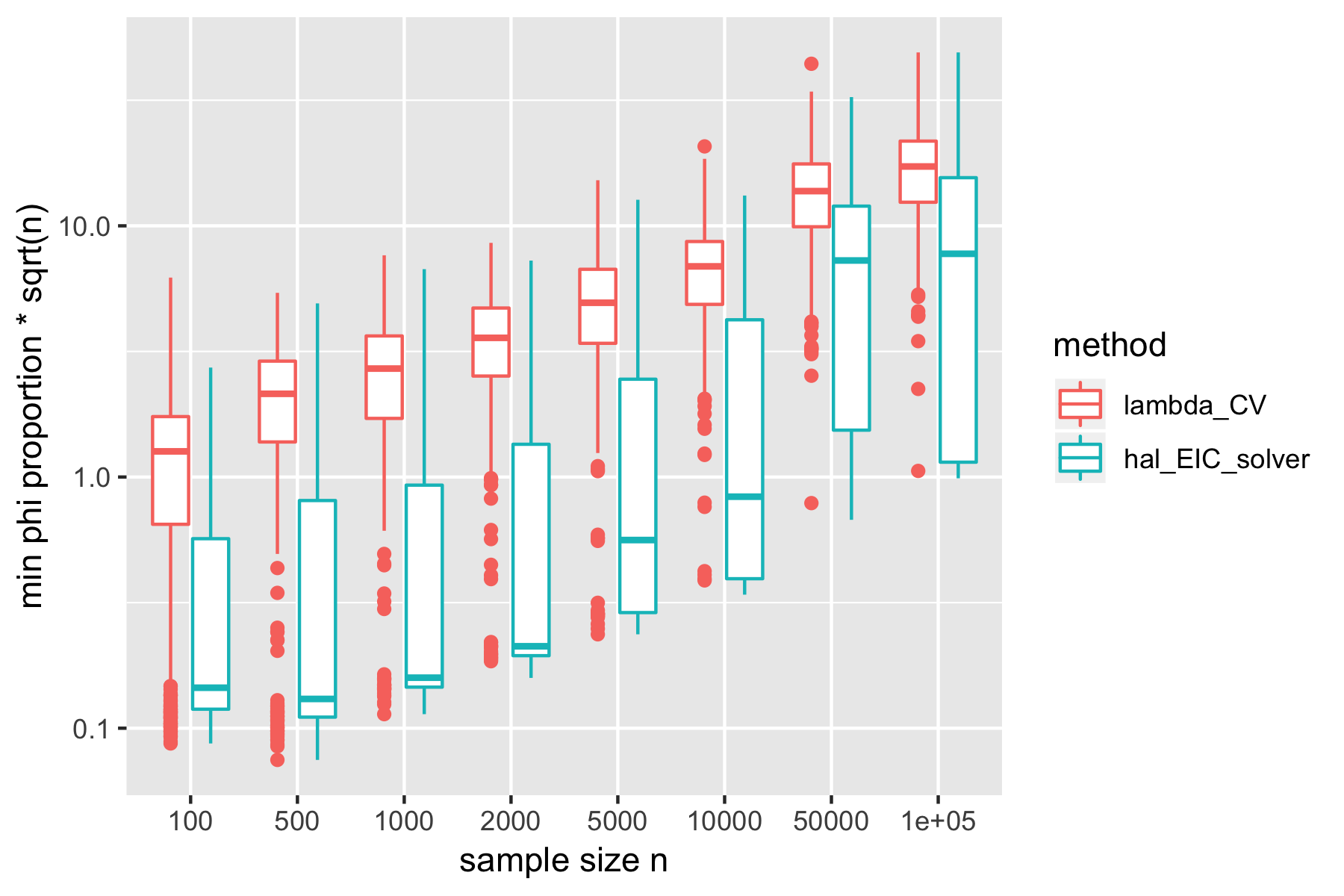}
    \caption{Summarization features for the average density value simulation: (a) $\sqrt{n}P_n(D^*)$ (b) $\sqrt{n}(P_n - P_0) (D_n - D_0)$ (c) $\|f_n\|$ (d) $\sqrt{n}\min_{s,j \in J_n} P_n \phi_{s,j}$ for undersmoothed HAL and HAL(CV) at different sample sizes}
    \label{fig:avgdens2}
\end{figure}

\section{Discussion}\label{section7}
In this article we established that for realistic and nonparametric statistical models an overfitted Spline HAL-MLE of a functional parameter of the data distribution results in efficient plug-in estimators of pathwise differentiable functionals of this functional parameter. The statistical model can be any model for which the parameter space for the functional is a (cartesian product of a) subset of the the set of multivariate cadlag functions with a universal bound on the sectional variation norm. 
The undersmoothing condition requires that one chooses the $L_1$-norm in the HAL-MLE large enough so that the basis functions with non-zero coefficients includes ''sparse enough'' basis functions, where ''sparse enough'' corresponds with assuming that the proportion of non-zero elements (among $n$ observations of this basis function) in the basis function converges to zero at a rate faster than a rate slower than $n^{-1/2}$. This rate could be as slow as $n^{-1/4}$ if one would be able to establish that the HAL-MLE converges in supremum norm at a rate faster than $n^{-1/4}$ as it does in $L^2(P_0)$-norm or loss-based dissimilarity, but, either way, the rate could be set at level $n^{-1/2+\alpha(k_1)}$, where $n^{-1/2-\alpha(k_1)}$ is the rate of the HAL-MLE w.r.t. loss-based dissimilarity. The undersmoothing condition represents a rate that is not parameter specific, so that such an undersmoothed HAL-MLE will be efficient for any of its smooth functionals. In addition, the undersmoothing of the HAL-MLE does not change its rate of convergence w.r.t. the HAL-MLE optimally tuned with cross-validation, suggesting that it is still a good estimator of the true functional parameter. 

On the other hand, a  typical TMLE targeting one particular target parameter will generally only be asymptotically efficient for that particular target parameter, and not even asymtotically linear for other smooth functionals, even if it uses as initial estimator the HAL-MLE tuned with cross-validation. Therefore it appears to be an interesting topic to better understand the sampling distribution of the undersmoothed HAL-MLE in an asymptotic sense and in relation to a sampling distribution of a TMLE using an optimally smoothed (i.e., cross-validation) HAL-MLE as initial estimator. Note, however, that if the TMLE uses an undersmoothed HAL-MLE as initial estimator, than the TMLE step should result in small changes, thereby mostly preserving the behavior of the undersmoothed HAL-MLE. 

It is also of interest to observe that the second order remainder of the HAL-MLE for a pathwise differentiable functional appears to either be driven by the square of the $L^2(P_0)$-norm of the HAL-MLE itself w.r.t. the functional parameter, or, in the case that the efficient influence curve has a nuisance parameter $G$, a second order remainder might also (or only) involve a product of differences of the HAL-MLE $Q_n$ w.r.t. its true counterpart $Q_0$ and the difference of a projection $G_{0,n}$ of the true $G_0$ w.r.t. the linear space of basis functions selected by the undersmoothed  HAL-MLE $Q_n$.
Since $G_{0,n}$ is a type of oracle estimator of $G_0$, this suggest that in a model in which our knowledge on $G_0$ is not any better than our knowledge on $Q_0$, this HAL-MLE has a good second order remainder that might generally be smaller than what it would be for a TMLE that  estimates $G_0$ with an actual estimator such as the HAL-MLE.

On the other hand, if the statistical model involves particularly strong knowledge on the nuisance parameter $G_0$, then a TMLE can fully utilize this model on $G_0$ and thereby obtain a better behaved second order remainder than the one for the overfitted HAL-MLE.
One also suspects that a TMLE will be more sensitive to lack of positivity for the target parameter than the undersmoothed HAL-MLE. 
Therefore, we conjecture that an undersmoothed HAL-MLE might be the preferred estimator in models in which case the estimation of $G_0$ is as hard as estimation of $Q_0$, and when lack of positivity is a serious issue,  while an HAL-TMLE might be the preferred estimator when estimation of $G_0$ is easier than estimation of $Q_0$. These are not formal statements, but indicate a qualitative  comparison between the undersmoothed HAL-MLE and a HAL-TMLE using an estimator (HAL-MLE) $G_n$ of $G_0$.

In future research we hope to address the comparison between undersmoothed HAL-MLE and HAL-TMLE in realistic simulations and possibly by formal comparison by their second order remainders. In a subsequent article we will marry the TMLE with the HAL-MLE by defining a targeted HAL-MLE that minimizes the empirical risk over the linear span of basis functions (approximating the true cadlag function with finite sectional variation norm) under the $L_1$-constraint {\em and} under the constraint that the Euclidean norm of the empirical mean of the efficient influence curve at the HAL-MLE (as well as at an estimator $G_n$) is $o_P(n^{-1/2})$.
We will show that undersmoothing this targeted HAL-MLE results in an estimator that is still efficient across all smooth functionals, while it is able to fully exploit all knowledge on $G_0$ for the sake of the specific target parameter. 

A key advantage of a TMLE is that it can utilize any super-learner so that its library can include many other powerful machine learning algorithms beyond $m$-th order Spline-HAL MLEs. 
In this manner a TMLE using a powerful super-learner might compensate for the favorable property of an undersmoothed HAL-MLE w.r.t. size of the second order remainder. In another future article we will provide a method that marries a powerful super-learner with HAL-MLE, by using the super-learner as a dimension reduction, and applying HAL-MLE as the meta learning step in an ensemble learner. We will show that an undersmoothed HAL-MLE in this metalearning step will result again in an estimator that is efficient for any of its smooth functionals. By actually using a targeted HAL-MLE as meta learning step, we might end up with an estimator that is able to still fully exploit the strengths super-learning, undersrmoothed HAL-MLE, and TMLE using a good esitmator of $G_0$, combined in one method.  

\section*{Acknowledgments}
Research reported in this publication was supported by the National Institute Of Allergy And Infectious Diseases of the National Institutes of Health under Award Number R01AI074345. The content is solely the responsibility of the authors and does not necessarily represent the official views of the National Institutes of Health.

\bibliography{combined,TLB2,CCbib,TMLEivan,gencens2003,gencens2003logfluct,TMLE1_TLB,TMLE_TLB,gencens2003a,gencens2003b,gencens2003c,gencens2003d,gencens2003e,Grant}

\appendix 

\section*{Appendix}

\section{Representing a function as an infinitesimal linear combination of spline-basis functions}
We have the following representation theorem for the smoothness class $D^m[0,\tau]$ consisting of $k$-variate real valued cadlag functions for which the $m$-th order sectional variation norm is bounded (and thereby is, in particular, $m$-times differentiable), as defined in Section 2. 

\begin{theorem}{\bf  $m$-th order spline representation of a function $Q\in D^m[0,\tau]$.}\label{thsplinerepr}
For any function $Q\in D^m[0,\tau]$ (i.e., finite $m$-th order sectional variation norm), we have
\begin{eqnarray*}
Q(x)&=&Q(0)+\sum_{j=0}^{m-1}\sum_{\bar{s}(j)}Q_{\bar{s}(j)}^{j+1}(0_{s_j})\phi_{\bar{s}(j),\emptyset,x_s}(0)\\
&&+\sum_{\bar{s}(m)}\int \phi_{\bar{s}(m),x_s}(z_{s_m})dQ^m_{\bar{s}(m)}(z_{s_m}).
\end{eqnarray*}

Thus, if $Q\in D^{m+1}[0,\tau]$, then we have
\begin{eqnarray*}
Q(x)&=&Q(0)+\sum_{j=0}^{m-1}\sum_{\bar{s}(j)}Q_{\bar{s}(j)}^{j+1}(0_{s_j})\phi_{\bar{s}(j),\emptyset,x_s}(0)\\
&&\hspace*{-1cm}
+\sum_{\bar{s}(m)} Q_{\bar{s}(m)}^{m+1}(0_{s_m})\phi_{\bar{s}(m),\emptyset,x_s}(0)
+\sum_{\bar{s}(m+1)}\int \phi_{\bar{s}(m+1),x_s}(z_{s_{m+1}}) dQ^{m+1}_{\bar{s}(m+1)}(z_{s_{m+1}}).
\end{eqnarray*}
\end{theorem}
{\bf Proof of Theorem:}

We already expressed $Q$ in terms of it integrals w.r.t. the measures generated by its $s$-specific sections $Q_s(x)=Q(x_s,0_{-s})$ for $s\subset\{1,\ldots,k\}$. Suppose that for each subset $s$ $Q_s$ is absolutely continuous w.r.t. Lebesque measure so that we have $dQ_s(u_s)=Q_s^1(u_s)du_s$. Suppose now that $Q_s^1$ is a cadlag function so that we can apply the same representation to $Q_s^1$:
\[
Q_s^1(u_s)=Q_s^1(0_s)+\sum_{s_1\subset s}\int_{(0_{s_1},u_{s_1}]} dQ_{s,s_1}^1(y_{s_1}).\]
In this manner we obtain a representation for $Q(x)$ in terms of integrals w.r.t. $Q_{s,s_1}$ across all subsets $s,s_1$ with $s_1\subset s$.
Specifically, we have
\begin{eqnarray*}
Q(x)&=&Q(0)+\sum_{s\subset\{1,\ldots,d\}}\int_{(0_s,x_s]}Q_s^1(u_s)du_s\\
&=&Q(0)+\sum_s\int_{(0_s,x_s]}Q_s^1(0_s)du_s+\sum_s\int_{(0_s,x_s]}\sum_{s_1\subset s}\int_{(0_{s_1},u_{s_1}]} dQ_{s,s_1}^1(y_{s_1}) du_s\\
&=&Q(0)+\sum_sQ_s^1(0_s)\prod_{j\in s}x_j \\
&&+\sum_{s,s_1\subset s}\int_{y_{s_1}}\left\{ \int_{u_s} I(u_s\leq x_s)I(y_{s_1}\leq u_{s_1}) du_s\right\}  dQ_{s,s_1}^1(y_{s_1}),
\end{eqnarray*}
where we used Fubini's theorem to exchange the order of integration. The product of the indicators can be written as 
$I(y_{s_1}\leq u_{s_1}\leq x_{s_1})I(y_{s_1}\leq x_{s_1})\prod_{j\in s,j\not\in s_1}I(u_j\leq x_j)$.
The inner integral represents a mapping from the original indicator basis function $x_s\rightarrow I(u_s\leq x_s)$ to a new basis function
\begin{eqnarray*}
x_s&\rightarrow& \phi_{s,s_1,x_s}(y_{s_1})\\
&\equiv& \int_{u_s} I(u_s\leq x_s)I(y_{s_1}\leq u_{s_1}) du_s=\prod_{j\in s_1}(x(j)-y(j))I(y(j)\leq x(j))\prod _{j\in s/s_1}x_j.\end{eqnarray*}
This is a  tensor product of spline basis functions across the components in $s$. Specifically, for $j\in s_1$ it is the first order spline-basis (line with slope 1 starting at knot $y_j$), while for 
$j\in s/s_1$, it is first order spline at knotpoint $0$.
We also define
\[
\phi_{s,\emptyset,x_s}(0_s)=\prod_{j\in s}x_j,\]
which corresponds with setting $s_1$ equal to empty set in definition of $\phi_{s,s_1,x_s}(y_s)$ and selecting knot $y_s=0_s$.
Thus, we have obtained:
\begin{eqnarray*}
Q(x)&=&Q(0)+\sum_sQ_s^1(0_s)\phi_{s,\emptyset,x_s}(0_s)+\sum_{s,s_1}\int_{y_{s_1}}\phi_{s,s_1,x_s}(y_{s_1}) dQ_{s,s_1}^1(y_{s_1}).\end{eqnarray*}
This proves the representation for functions in $D^m[0,\tau]$ for $m=1$.
This representation shows that this class of functions can be represented as an infinitesimal linear combination of first order spline basis functions for which the $L_1$-norm of the coefficients equals the first-order sectional variation norm of the function.


Let's now use the same method to derive a representation by assuming another degree of differentiability, thereby establishing the general story. The last expression expresses $Q$ in integrals w.r.t $Q_{s,s_1}^1$. Suppose now that $Q_{s,s_1}^1$ is absolutely continuous w.r.t. Lebesque measure so that $dQ_{s,s_1}^1(y_{s_1})=Q_{s,s_1}^2(y_{s_1})dy_{s_1}$ for a second order  derivative $Q_{s,s_1}^2$. Let's now also assume that $Q_{s,s_1}^2 $ is cadlag and has finite sectional variation norm so that\[Q_{s,s_1}^2(y_{s_1})=Q_{s,s_1}^2(0_{s_1})+\sum_{s_2\subset s_1}\int_{(0_{s_2},y_{s_2}]} dQ_{s,s_1,s_2}^2(z_{s_2}).\]Substitution of this into the last expression for $Q$, combined with the derivation above involving change or order of integration, yields the following:
\begin{eqnarray*}
Q(x)&=&Q(0)+\sum_sQ_s^1(0_s)\phi_{s,\emptyset,x_s}(0_s)+
\sum_{s,s_1}Q_{s,s_1}^2(0_{s_1})\phi_{s,s_1,\emptyset,x_s}(0)\\
&&
+
\sum_{s,s_1,s_2}\int_{z_{s_2}} \phi_{s,s_1,s_2,x_s}(z_{s_2}) dQ_{s,s_1,s_2}^2(z_{s_2}).
\end{eqnarray*}
The above representation shows that this function is represented as an infinitesimal linear combination of (up till) second order spline basis functions for which the $L_1$-norm of the coefficients  equals the second order sectional variation norm of the function. Some of the spline basis functions in the tensor product basis functions are of lower order, but only at knot points equal to $0$. 

Similarly, the 3rd  order representation follows from the above 2-nd order representation:
\begin{eqnarray*}
Q(x)&=&Q(0)+\sum_sQ_s^1(0_s)\phi_{s,\emptyset,x_s}(0_s)+
\sum_{s,s_1}Q_{s,s_1}^2(0_{s_1})\phi_{s,s_1,\emptyset,x_s}(0)\\
&&+
\sum_{s,s_1,s_2}Q_{s,s_1,s_2}^3(0_{s_2})\phi_{s,s_1,s_2,\emptyset,x_s}(0)\\
&&+
\sum_{s,s_1,s_2,s_3}\int_{z_{s_3}} \phi_{s,s_1,s_2,s_3,x_s}(z_{s_3}) dQ_{s,s_1,s_2,s_3}^3(z_{s_3}),
\end{eqnarray*}
where  the second order basis functions $x_s\rightarrow \phi_{s,s_1,s_2,x_s}(y_{s_2})$ indexed  by knots $y_2$  is mapped to the new third order basis functions
\[
\phi_{\bar{s}(3),x_s}(z_{s_3})\equiv \prod_{j\in s_3}\int_{(z_j,x_j]}\phi_{j,\bar{s}(2),x_j}(y_j)dy_j\prod_{j\in s_2/s_3}\int_{(0,x_j]}\phi_{j,\bar{s}(2),x_j}(y_j)dy_j
\prod_{j\in s/s_2}\phi_{j,\bar{s}(2),x_j}.
\]
In general, the $m$-th order spline representation of $Q$ is represented as follows:
\begin{eqnarray*}
Q(x)&=&Q(0)+\sum_{j=1}^{m-1}\sum_{\bar{s}(j)}Q_{\bar{s}(j)}^{j+1}(0_{s_j})\phi_{\bar{s}(j),\emptyset,x_s}(0) +\sum_{\bar{s}(m)}\int\phi_{\bar{s}(m),x_s}(z_{s_m}) dQ^m_{\bar{s}(m)}(z_{s_m}),
\end{eqnarray*}
and the $m+1$-th order spline representation is derived from this one as follows:
\begin{eqnarray*}
Q(x)&=&Q(0)+\sum_{j=1}^{m-1}\sum_{\bar{s}(j)}Q_{\bar{s}(j)}^{j+1}(0_{s_j})\phi_{\bar{s}(j),\emptyset,x_s}(0)\\
&&+\sum_{\bar{s}(m)} Q_{\bar{s}(m)}^{m+1}(0_{s_m})\phi_{\bar{s}(m),\emptyset,x_s}(0)
+\sum_{\bar{s}(m+1)}\int \phi_{\bar{s}(m+1),x_s}(z_{s_{m+1}}) dQ^{m+1}_{\bar{s}(m+1)}(z_{s_{m+1}}),
\end{eqnarray*}
where  the $m$-th order basis functions $x_s\rightarrow \phi_{\bar{s}(m),x_s}(y_{s_m})$  indexed by knots $y_{s_m}$ is mapped to the new $m+1$-th order basis functions
\begin{eqnarray*}
\phi_{\bar{s}(m+1),x_s}(z_{s_{m+1}})
&\equiv& \prod_{j\in s_{m+1}}\int_{(z_j,x_j]}\phi_{j,\bar{s}(m),x_j}(y_j)dy_j\prod_{j\in s_{m}/s_{m+1}}\int_{(0,x_j]}\phi_{j,\bar{s}(m),x_j}(y_j)dy_j\\
&&
\prod_{j\in s/s_m}\phi_{j,\bar{s}(m),x_j}(0).
\end{eqnarray*}
Note that the last term in the $m$-th order representation is replaced by the last two new terms to obtain the $m+1$-th order representation. This completes the proof.
$\Box$

\section{Rate of convergence of $m$-th order Spline HAL-MLE, and smoothness-adaptive Spline HAL-MLE}
Let ${\cal Q}^m=D^m[0,\tau]$.
Under the assumption that $M_1<\infty$ and $M_{20}<\infty$, we have that
$r_n(m)\equiv d_0(Q_n^m,Q_0^m)=o_P(n^{-1/2})$, where $Q_0^m=\arg\min_{Q\in {\cal Q}^m}P_0L(Q)$. Let $m_0$ be the smallest integer $m$ for which $Q_0\in {\cal Q}^m$. Then, we have that, if $m\geq m_0$, $d_0(Q_n^m,Q_0)=r_n(m)=o_P(n^{-1/2})$. In general, the rate of convergence $r_n(m)$ will be unique for each $m$ so that the highest rate of convergence is achieved by the $m_0$-th order Spline HAL-MLE $Q_n^{m_0}$, which is achieved by selecting $m$ with cross-validation (due to asymptotic equivalence of cross-validation selector with oracle selector). 
Thus, if $m_n=\arg\min_m E_{B_n}P_{n,B_n}^1 L(\hat{Q}^m(P_{n,B_n}^0))$ is the cross-validation selector, then, by the oracle inequality of the cross-validation selector, we have
 \begin{equation}\label{asympequiv}
\frac{ E_{B_n}d_0(\hat{Q}_n^{m_n}(P_{n,B_n}^0),Q_0)}{\min_mE_{B_n}d_0(\hat{Q}^{m}(P_{n,B_n}^0),Q_0)}\rightarrow_p 1.\end{equation}
 We refer to $Q_n^{m_n}=\hat{Q}^{m_n}(P_n)$ as the smoothness adaptive Spline HAL-MLE.

\subsection{Asymptotic equivalence of the smoothness adaptive Spline HAL-MLE with the oracle Spline HAL-MLE}
The following lemma states that if the rates of convergence of the $m$-th order Spline HAL-MLE are unique across $m$, and the conditions under which the cross-validation selector $m_n$ is asymptotically equivalent with the oracle selector (\ref{asympequiv}) hold, then it follows that $P(m_n=m_0)\rightarrow 1$. 
In that case, establishing that the plug-in $m_0$-th order Spline HAL-MLE is asymptotically efficient also implies that the plug in of the smoothness adaptive Spline HAL-MLE is asymptotically efficient. Thus, under this condition, our asymptotic efficiency results for the $m$-th order Spline HAL-MLE implies also our desired result for the asymptotic efficiency of the smoothness adaptive Spline HAL-MLE.

\begin{lemma}\label{equivm0}
Let $r_m(n)\equiv E_{B_n}d_0(\hat{Q}^m(P_{n,B_n}^0),Q^m_0)$.  Let $m_n\in \{1,\ldots,K_n\}$. Assume that the $m$-specific rates $r_m(n)$ are unique in the sense that
\begin{equation}\label{uniquerateassumption}
 \frac{r_{m_0}(n)}{\min_{m<m_0}r_m(n)}\rightarrow_p 0;\end{equation} 
  \begin{equation}\label{equiva}
 \frac{\log K_n/n}{\min_m r_m(n)}\rightarrow_p 0;
 \end{equation}
 and
 $M_1<\infty$, $M_{20}<\infty$.
 
 Then, (\ref{asympequiv}) holds, and
  $P(m_n=m_0)\rightarrow 1$ as $n\rightarrow\infty$. 
 \end{lemma}
{\bf Proof of Lemma:}
Due to the asymptotic equivalence of the cross-validation selector $m_n$ with the oracle selector under the above conditions, it follows that  $\lim  P(m_n>m_0)\rightarrow 0$. Suppose now that $\lim\sup_m P(m_n<m_0)>\delta>0$ for some $\delta>0$. Let $A=\{m:m>m_0\}$ so that $P(m_n\in A)>\delta>0$. Then, for $m_n\in A$, we have \[
r_{m_0}(n)/r_{m_n}(n)\leq r_{m_0}(n)/\min_{m>m_0}r_m(n),\] but  the latter upper-bound   converges to zero in probability, by assumption (\ref{uniquerateassumption}). Thus, in that case $r_{m_0}(n)/r_{m_n}(n)$ does not converge to $1$ in probability.
But $\min_m r_m(n)/r_{m_n}(n)\leq r_{m_0}(n)/r_{m_n}(n)$, so that this also implies that $\min_m r_m(n)/r_{m_n}(n)$ does not converge to $1$ in probability.
However, this contradicts (\ref{asympequiv}) which states that $\min_m r_m(n)/r_{m_n}(n)\rightarrow_p 1$.
Thus, we can conclude that $P(m_n=m_0)\rightarrow 1$ as $n\rightarrow\infty$. $\Box$

\section{Proof of Theorem \ref{thscoreequation} and Lemma \ref{lemmaassumptiona}.}
The HAL-MLE has the form $Q_n=\sum_{s,j\in {\cal J}_n(s)}\beta_n(s,j)\phi_{s,j}$ for  a finite collection of basis functions. A basis function $\phi_{s,j}(X) $ has support on 
$\{X:X_s>x_{s,j}\}$ for a knot point $x_{s,j}$, and $s$ a subset of $\{1,\ldots,k\}$.
We also know that $\sum_{s,j}\mid \beta_n(s,j)\mid  \leq C_n$ for the selected $L_1$-bound $C_n$ (typically  the $L^1$-norm will be equal to $C_n$).
We have that
\[
\beta_n=\arg\min_{\beta,\sum_{s,j}\mid \beta(s,j)\mid\leq C_n}P_n L\left(\sum_{s,j\in {\cal J}_n(s)}\beta(s,j)\phi_{s,j}\right).\]
Consider paths $(1+\epsilon h(s,j))\beta_n(s,j)$ for a  bounded vector $h$, which yields a collection of scores
\[
S_h(Q_n)=\frac{d}{dQ_n}L(Q_n)\left (\sum_{s,j\in {\cal J}_n(s)}h(s,j)\beta_n(s,j)\phi_{s,j}\right ).\]
Let $r(h,Q_n)=\sum_{s,j\in {\cal J}_n(s)} h(s,j)\mid \beta_n(s,j)\mid$. If $r(h,Q_n)=0$, then for $\epsilon$ small enough, 
\begin{eqnarray*}
\sum_{s,j\in {\cal J}_n(s)}\mid (1+\epsilon h(s,j))\beta_n(s,j)\mid&=&\sum_{s,j\in {\cal J}_n(s)}(1+\epsilon h(s,j))\mid \beta_n(s,j)\mid \\
&=&\sum_{s,j\in {\cal J}_n(s)}\mid \beta_n(s,j)\mid+\epsilon r(h,Q_n)\\
&=&\sum_{s,j\in {\cal J}_n(s)}\mid \beta_n(s,j)\mid .\end{eqnarray*}
Thus, by $\beta_n$ being an MLE,  $P_n S_h(Q_n)=0$ for any $h$ satisfying $r(h,Q_n)=0$.
Let $h^*=h^*_n$ be chosen so that $P_n S_{h^*_n}(Q_n)=P_n D^*_n(Q_n,G_0)$ for the  approximation $D^*_n(Q_n,G_0)$ of $D^*(Q_n,G_0)$ specified in the theorem.
We want to show that $P_n D^*_n(Q_n,G_0)=o_P(n^{-1/2})$, i.e. $P_n S_{h^*_n}(Q_n)=o_P(n^{-1/2})$.
Let $(s^*,j^*)$  be a particular choice in our finite index set ${\cal J}_n$ satisfying $\beta_n(s^*,j^*)\not =0$, which we can specify later to minimize the bound. 
Let $\tilde{h}$ be defined by $\tilde{h}(s,j)=h^*(s,j)$ for $(s,j)\not =(s^*,j^*)$, and $\tilde{h}(s^*,j^*)$ is defined by $r(\tilde{h},Q_n)=\sum_{s,j\in {\cal J}_n(s)} \tilde{h}(s,j)\mid \beta_n(s,j)\mid=0$, so that we know $P_n S_{\tilde{h}}(Q_n)=0$.
Thus,\[
\sum_{(s,j)\not =(s^*,j^*)}h^*(s,j)\mid \beta_n(s,j)\mid+ 
\tilde{h}(s^*,j^*)\mid \beta_n(s^*,j^*)\mid =0.\]
This gives
\[
\tilde{h}(s^*,j^*)=-\frac{\sum_{(s,j)\not =(s^*,j^*)}h^*(s,j)\mid \beta_n(s,j)\mid}{
\mid \beta_n(s^*,j^*)\mid }.\]
So
\[
\begin{array}{l}
\sum_{s,j}(\tilde{h}-h^*)(s,j)\beta_n(s,j)\phi_{s,j}=
(\tilde{h}-h^*)(s^*,j^*)\beta_n(s^*,j^*)\phi_{s^*,j^*}\\
=\left(-\frac{\sum_{(s,j)\not =(s^*,j^*)}h^*(s,j)\mid \beta_n(s,j)\mid}{
\mid \beta_n(s^*,j^*)\mid } \beta_n(s^*,j^*)- h^*(s^*,j^*)\beta_n(s^*,j^*) \right)\phi_{s^*,j^*}\\
\equiv c_n(s^*,j^*)\phi_{s^*,j^*},
\end{array}
\]
where
\[
c_n(s^*,j^*)=-\frac{\sum_{(s,j)\not =(s^*,j^*)}h^*(s,j)\mid \beta_n(s,j)\mid}{
\mid \beta_n(s^*,j^*)\mid}  \beta_n(s^*,j^*)- h^*(s^*,j^*)\beta_n(s^*,j^*).\]
We note that $c_n(s^*,j^*)$ is bounded by $\sum_{s,j} \mid h^*(s,j)\mid \mid \beta_n(s,j)\mid  $. So we can bound this by 
$\pl h^*\pl_{\infty} C_n$.
Thus under the assumption that $\pl h^*_n\pl_{\infty}=O_P(1)$, we have that $c_n(s^*,j^*)=O_P(1)$.

For this choice $\tilde{h}$, let's compute $P_n S_{\tilde{h}}(Q_n)-P_n S_{h^*}(Q_n)$ (which equals $P_n S_{h^*}(Q_n)$):
\begin{eqnarray*}
P_n S_{\tilde{h}}(Q_n)-P_n S_{h^*}(Q_n)&=&
P_n \frac{d}{dQ_n}L(Q_n)\left(\sum_{s,j}(\tilde{h}-h^*)(s,j)\beta_n(s,j)\phi_{s,j}\right)\\
&=&P_n \frac{d}{dQ_n}L(Q_n)( c_n(s^*.j^*)\phi_{s^*,j^*})\\
&=&c_n(s^*,j^*) P_n \frac{d}{dQ_n}L(Q_n)(\phi_{s^*,j^*})\\
&=&O_P\left( P_n \frac{d}{dQ_n}L(Q_n)(\phi_{s^*,j^*})\right).
\end{eqnarray*}
Therefore, our undersmoothing condition is that
\begin{equation}\label{keycondition}
\min_{s,j\in {\cal J}_n(s),\beta_n(s,j)\not =0}\pl P_n \frac{d}{dQ_n}L(Q_n)(\phi_{s,j})\pl =o_P(n^{-1/2}).\end{equation}
Under this condition we have $P_n S_{\tilde{h}}(Q_n)-P_n D^*_n(Q_n,G_0)=o_P(n^{-1/2})$, but, since $P_n S_{\tilde{h}}(Q_n)=0$, this implies the desired conclusion $P_n D^*_n(Q_n,G_0)=o_P(n^{-1/2})$.
This proves the first statement of Theorem \ref{thscoreequation}.

Let now $(s^*,j^*)=\arg\min_{s,j\in {\cal J}_n(s)}P_0 \phi_{s,j}$.
To understand, $P_n \frac{d}{dQ_n}L(Q_n)(\phi_{s^*,j^*})$ we can proceed as follows.
\begin{eqnarray*}
P_n \frac{d}{dQ_n}L(Q_n)(\phi_{s^*,j^*})&=&
(P_n-P_0)\frac{d}{dQ_n}L(Q_n)(\phi_{s^*,j^*})+P_0 \frac{d}{dQ_n}L(Q_n)(\phi_{s^*,j^*}).
\end{eqnarray*}
Let $S_{s,j}(Q_n)\equiv \frac{d}{dQ_n}L(Q_n)(\phi_{s,j})$.
Suppose that  $P_0 S_{s^*,j^*}(Q_n)^2\rightarrow_p 0$, which will generally hold whenever $P_0\phi_{s^*,j^*}=o_P(1)$. 
We also have that $\{S_{s,j}(Q):Q\in {\cal Q},(s,j)\}$ is contained in the class of cadlag functions with uniformly bounded sectional variation norm, which is a Donsker class.
Thereby, by asymptotic equicontinuity of the empirical process indexed by a Donsker class, we have $(P_n-P_0)S_{s^*,j^*}(Q_n)=o_P(n^{-1/2})$.
Thus, it remains to show that $P_0 S_{s^*,j^*}(Q_n)=o_P(n^{-1/2})$.
We now note that
\[
P_0 S_{s^*,j^*}(Q_n)=P_0 \{S_{s^*,j^*}(Q_n)-S_{s^*,j^*}(Q_0)\}+
P_0 S_{s^*,j^*}(Q_0),\]
but $P_0  S_{s,j}(Q_0)=0$ for all $(s,j)$, since $Q_0=\arg\min_{Q}P_0L(Q)$.
Therefore, $P_n \frac{d}{dQ_n}L(Q_n)(\phi_{s^*,j^*})=o_P(n^{-1/2})$ if
\begin{equation}\label{suffcritical}
P_0\{S_{s^*,j^*}(Q_n)-S_{s^*,j^*}(Q_0)\}=o_P(n^{-1/2}).\end{equation}
This proves the second statement of Theorem \ref{thscoreequation}.
The third statement is a trivial implication, which completes the proof of Theorem \ref{thscoreequation}. $\Box$

{\bf Proof of Lemma \ref{lemmaassumptiona}:}
Consider the special case that $O=(Z,X)$, $L(Q)(O)$ depends on $Q$ through $Q(X)$ only, and $\frac{d}{dQ}L(Q)(\phi)=\frac{d}{dQ}L(Q)\times \phi$, i.e., the directional derivative 
$\left . \frac{d}{d\epsilon}L(Q+\epsilon \phi)\right |_{\epsilon =0}$ of $L()$ at $Q$ in the direction $\phi$ is just multiplication of a function $\frac{d}{dQ}L(Q)$ of $O$ with 
$\phi(X)$. In that case, we have that (\ref{suffcritical}) reduces to
\begin{equation}\label{suffcriticala}
P_0\left\{\frac{d}{dQ_n}L(Q_n)-\frac{d}{dQ_0}L(Q_0) \right\} \phi_{s^*,j^*}=o_P(n^{-1/2}).\end{equation}
We assume
\[
\pl \frac{d}{dQ_n}L(Q_n)-\frac{d}{dQ_0}L(Q_0)\pl_{\infty}=O(\pl Q_n-Q_0\pl_{\infty}).\]
Then, (\ref{suffcriticala}) reduces to
\[
\pl Q_n-Q_0\pl_{\infty}P_0 \phi_{s^*,j^*} =o_P(n^{-1/2}).\]
This teaches us  that the critical condition (\ref{keycondition}) holds if \[
\min_{s,j\in {\cal J}_n(s),\beta_n(s,j)\not =0}P_0\phi_{s,j}=O_P(n^{-1/2}). \]
and that for this choice $(s^*,j^*)$ we have $P_0\{S_{s^*,j^*}(Q_n)\}^2\rightarrow_p 0$.
The latter holds if $\min_{s,j}P_0\phi_{s,j}=o_P(1)$, since $\frac{d}{dQ_n}L(Q_n) $ is uniformly bounded. Finally, since
 $P_0\phi_{s,j}=O(P_0(X(s)\geq x_{s,j}))$, $\sup_{s,j}\mid (P_n-P_))\phi_{s,j}\mid=O_P(n^{-1/2})$, we can replace $P_0\phi_{s,j}$ by $P_n \phi_{s,j}$ in the condition. This proves Lemma \ref{lemmaassumptiona}.

$\Box$

 \section{Proof of Theorem \ref{thfinalspline}}

Let $G_{0n}$ be an approximation of $G_0$, and let $D^*(Q_n,G_{0n})$ be the approximation of $D^*(Q_n,G_0)$ in the space of scores ${\cal S}(Q_n)$.
 We have the following general theorem which proves the first part of Theorem \ref{thfinalspline}. \begin{theorem}\label{thgen}
 Consider the  HAL-MLE $Q_n$ with $C=C_u$ or $C=C_n$.  Assume $M_1,M_{20}<\infty$. We have $d_0(Q_n,Q_0)=O_P(n^{-1/2-\alpha(k_1)})$.
  Assume also that for a given approximation $G_{0n}\in {\cal G}$ of $G_0$ which satisfies
  \begin{equation}\label{effscore}
  P_n D^*(Q_n,G_{0n})=o_P(n^{-1/2}).
  \end{equation}
  \begin{itemize}
  \item $R_2((Q_n,G_{0n}),(Q_0,G_0))=o_P(n^{-1/2})$ and 
$P_0\{D^*(Q_n,G_{0n})-D^*(Q_0,G_0)\}^2\rightarrow_p 0$.
  \item  $\{D^*(Q,G):Q\in {\cal Q},G\in {\cal G}\}$ is contained in the class of $k_1$-variate cadlag functions on a cube $[0,\tau_o ]\subset\openr^{k_1}$ in a Euclidean space and that 
$\sup_{Q\in {\cal Q},G\in {\cal G}}\pl D^*(Q,G)\pl_v^*<\infty$.  
\end{itemize}
Then $\Psi(Q_n)$ is asymptotically efficient at $P_0$.
\end{theorem}
{\bf Proof:}
The exact second order expansion at $G_{0n}$ of the target parameter $\Psi$ yields
\begin{eqnarray*}
\Psi(Q_n)-\Psi(Q_0)&=&(P_n-P_0)D^*(Q_n,G_{0n})-P_n D^*(Q_n,G_{0n})\\
&&\hfill + R_2((Q_n,G_{0n}),(Q_0,G_0)).\end{eqnarray*}
Given that $d_0(Q_n,Q_0)=O_P(n^{-1/2-\alpha(k_1)})$, and that $G_{0n} $ is presumably at least as good of an approximation of $G_0$, it is a reasonable assumption to assume  $R_2((Q_n,G_{0n}),(Q_0,G_0))=o_P(n^{-1/2})$ and 
$P_0\{D^*(Q_n,G_{0n})-D^*(Q_0,G_0)\}^2\rightarrow_p 0$.
We also assume that $\{D^*(Q,G):Q\in {\cal Q},G\in {\cal G}\}$ is contained in the class of $k_1$-variate cadlag functions on a cube $[0,\tau_o ]\subset\openr^{k_1}$ in a Euclidean space and that 
$\sup_{Q\in {\cal Q}, G\in {\cal G}}\pl D^*(Q,G)\pl_v^*<\infty$. 
This essentially states that the sectional variation norm of $D^*(Q,G)$ can be bounded in terms of the sectional  variation norm of $Q$ and $G$, which will naturally hold under a strong positivity assumption that bounds denominators away from zero.
Since the class of cadlag functions on $[0,\tau_o]$ with sectional variation norm bounded by a universal constant is a Donsker class, 
 empirical process theory  yields:
 \[
 \Psi(Q_n)-\Psi(Q_0)=(P_n-P_0)D^*(Q_0,G_0) -P_n D^*(Q_n,G_{0n})+o_P(n^{-1/2}).\Box
 \]
 
This theorem can be easily generalized to a more general approximation $D^*_n(Q_n,G_0)\in {\cal S}(Q_n)$  of $D^*(Q_n,G_0)$ (not necessarily of form  $D^*_n(Q_n,G_0)=D^*(Q_n,G_{0n})$ for some $G_{0n}$).

\begin{theorem}\label{thgena}
 Consider the  HAL-MLE $Q_n$ with $C=C_u$ or $C=C_n$.  Assume $M_1,M_{20}<\infty$. We have $d_0(Q_n,Q_0)=O_P(n^{-1/2-\alpha(k_1)})$.
  Assume also that for a given approximation $D^*_n(Q_n,G_0)$ we have
  $P_n D^*(Q_n,G_{0n})=o_P(n^{-1/2})$. In addition, assume
  \begin{itemize}
  \item $R_2((Q_n,G_{0}),(Q_0,G_0))=o_P(n^{-1/2})$, $P_0\{D^*_n(Q_n,G_0)-D^*(Q_n,G_0)\}=o_P(n^{-1/2})$, and 
$P_0\{D^*_n(Q_n,G_{0})-D^*(Q_0,G_0)\}^2\rightarrow_p 0$.
  \item  $\{D^*_n(Q,G_0),D^*(Q,G_0):Q\in {\cal Q}\}$ is contained in the class of $k_1$-variate cadlag functions on a cube $[0,\tau_o ]\subset\openr^{k_1}$ in a Euclidean space and that 
$\sup_{Q\in {\cal Q}}\max(\pl D^*(Q,G_0)\pl_v^*,\pl D^*_n(Q,G_0)\pl_v^*)<\infty$.  
\end{itemize}
Then $\Psi(Q_n)$ is asymptotically efficient at $P_0$.
\end{theorem}

  
Therefore, in order to prove Theorem \ref{thfinalspline}, it remains  to establish the condition under which  (\ref{effscore}) holds, which was proven in the previous Appendix.
 
\subsection{General proof of efficient score equation condition at $G_0$}
This subsection can be skipped for the purpose of proving Theorem \ref{thfinalspline}, but the following result fits here.

 \begin{lemma}\label{lemma}
  Under the conditions of Theorem \ref{thgen}, if $P_n D^*(Q_n,G_{0n})=o_P(n^{-1/2})$, then also $P_n D^*(Q_n,G_0)=o_P(n^{-1/2})$.
 Under the conditions of  Theorem \ref{thgena}, if $P_n D^*_n(Q_n,G_0)=o_P(n^{-1/2})$, then also $P_n D^*(Q_n,G_0)=o_P(n^{-1/2})$.
  \end{lemma}
 {\bf Proof:}
Firstly, we have
\begin{eqnarray*}
P_n D^*(Q_n,G_0)&=&P_n D^*_n(Q_n,G_{0})+P_n \{D^*(Q_n,G_0)-D^*_n(Q_n,G_{0})\}\\
&=&P_n \{D^*(Q_n,G_0)-D^*_n(Q_n,G_{0})\}+o_P(n^{-1/2}) .\end{eqnarray*}
In addition, we have
\begin{eqnarray*}
P_n \{D^*(Q_n,G_0)-D^*_n(Q_n,G_{0})\}&=&(P_n-P_0) \{D^*(Q_n,G_0)-D^*_n(Q_n,G_{0})\}\\
&&+P_0\{D^*(Q_n,G_0)-D^*_n(Q_n,G_{0})\}\\
&=& o_P(n^{-1/2})+ P_0\{D^*(Q_n,G_0)-D^*_n(Q_n,G_{0})\},
\end{eqnarray*}
since $\sup_{Q\in Q({\cal M})}\max(\pl D^*(Q,G_0)\pl_v^*,\pl D^*_n(Q,G_0)\pl_v^*) <\infty$, and $P_0 \{D^*_n(Q_n,G_{0})-D^*(Q_n,G_0)\}^2\rightarrow_p 0$.
If $D^*_n(Q_n,G_0)=D^*(Q_n,G_{0n})$, then the first assumption holds if $\sup_{P\in {\cal M}}\pl D^*(P)\pl_v^*<\infty$.

To understand the last term, consider the case that $D^*_n(Q_n,G_0)=D^*(Q_n,G_{0n})$. By the exact second order expansion $\Psi(Q_n)-\Psi(Q_0)=-P_0D^*(Q_n,G)+R_{20}(Q_n,G,Q_0,G_0)$ for all $G$, we have
\[
P_0\{D^*(Q_n,G_0)-D^*(Q_n,G_{0n})\}=R_{20}(Q_n,G_0,Q_0,G_0)-R_{20}(Q_n,G_{0n},Q_0,G_0).\]
In our general theorem \ref{thgen} we assumed $R_{20}(Q_n,G_{0n},Q_0,G_0)=o_P(n^{-1/2})$, which certainly implies $R_{20}(Q_n,G_0,Q_0,G_0)$ (which actually equals zero in double robust problems). 
 This then establishes that
 \[
 P_n D^*(Q_n,G_0)=o_P(n^{-1/2}).\]
 For general $D^*_n(Q_n,G_0)$, Theorem \ref{thgena} simply assumed $P_0 \{D^*(Q_n,G_0)-D^*_n(Q_n,G_0)\}=o_P(n^{-1/2})$.
 $\Box$

\section{Efficiency of HAL-MLE for general non-linear risk functions.}


{\bf Formulation of statistical estimation problem:}
Consider  a smooth parameter $\Psi:{\cal M}\rightarrow\openr^d$, where $\Psi(P)=\Psi(Q(P))$ for a $Q(P_0)$ that is identified as the  minimizer of a risk function 
$R(Q,P_0)$ over all $Q\in Q({\cal M})$: $Q_0=\arg\min_{Q\in Q({\cal M})}R(Q,P_0)$. As above, we assume that $\Psi$ is pathwise differentiable at $P$ with canonical gradient $D^*(Q(P),G(P))$ for a nuisance parameter $G(P)$, and, for each $Q\in Q({\cal M})$,  $P\rightarrow R_Q(P)\equiv R(Q,P)$ is pathwise differentiable at $P$ with canonical gradient $D^*_Q(P)$. Let $R_2(P,P_0)=\Psi(P)-\Psi(P_0)+P_0D^*(P)$, and
$R_{2Q}(P,P_0)=R_Q(P)-R_Q(P_0)+P_0D^*_Q(P)$ be the exact second order remainders for these two pathwise differentiable target parameters.

{\bf Defining the general HAL-MLE for non-linear risk functions:}
Let $Q_n=\arg\min_{Q\in Q({\cal M}),\pl Q\pl_v^*<C_n}R(Q,P_n^*)$ be the risk-based HAL-MLE defined above, where $P_n^*=P_{n,Q}^*$ is  an estimator of $R(Q,P_0)$ satisfying
$P_n D^*_Q(P_{n,Q}^*)=o_P(n^{-1/2})$, such as a TMLE (or an undersmoothed HAL-MLE).
Notice that we focus here on generalizing the $m=0$-Spline HAL-MLE. 

\subsection{$n^{-1/2}$-Consistency of the general HAL-MLE for non-linear risk functions.}
A generalization of the proof of the rate of convergence for the HAL-MLE w.r.t. loss-based dissimilarity establishes that $d_0(Q_n,Q_0)\equiv R(Q_n,P_0)-R(Q_0,P_0)=o_P(n^{-1/2})$.
\begin{lemma}
We will write $R(Q_n,P_n^*)$ for the plug in estimator of $R(Q_n,P_0)$ and $R(Q_n,P_{0,n}^*)$ for the plug-in estimator  of $R(Q_0,P_0)$ (treating $Q_0$ as given), satisfying $P_n D^*_{Q_n}(P_n^*)=o_P(n^{-1/2})$ and $P_n D^*_{Q_0}(P_{0,n}^*)=o_P(n^{-1/2})$, respectively.
(We  could have that $P_n^*$ and $P_{0,n}^*$ are TMLEs targeting $R(Q_n,P_0)$ and $R(Q_0,P_0)$, respectively, but we could also have that $P_n^*=P_{0,n}^*$ is an undersmoothed HAL-MLE that solves both efficient influence curve equations.)

We make the following assumptions:
\begin{itemize}
\item $R(Q_0,P_n^*)-R(Q_0,P_{0,n}^*)=o_P(n^{-1/2})$.
\item $\{D^*_Q(P):Q\in Q({\cal M}),P\in {\cal M}\}$ is a $P_0$-Donsker class.
\item $P_0\{D^*_{Q_n}(P_n^*)-D^*_{Q_0}(P_{0,n}^*)\}\rightarrow_p 0$, where we can use that 
$d_0(Q_n,Q_0)\rightarrow_p 0$.
\item $R_{2Q_n}(P_n^*,P_0)-R_{2Q_0}(P_{0,n}^*,P_0)=o_P(n^{-1/2})$.
\end{itemize}
Then, $d_0(Q_n,Q_0)=o_P(n^{-1/2})$.
\end{lemma}
{\bf Proof:}
Let $d_0(Q_n,Q_0)=R(Q_n,P_0)-R(Q_0,P_0)$. We have
\begin{eqnarray*}
0&\leq& d_0(Q_n,Q_0)\\
&=& R(Q_n,P_0)-R(Q_0,P_0)\\
&=&R(Q_n,P_0)-R(Q_n,P_n^*)-(R(Q_0,P_0)-R(Q_0,P_{0,n}^*))\\
&&+R(Q_n,P_n^*)-R(Q_0,P_{n}^*)\\
&&+R(Q_0,P_n^*)-R(Q_0,P_{0,n}^*)\\
&\leq& R(Q_n,P_0)-R(Q_n,P_n^*)-R(Q_0,P_0)+R(Q_0,P_{0,n}^*)\\
&&+R(Q_0,P_n^*)-R(Q_0,P_{0,n}^*)\\
&=&R(Q_n,P_0)-R(Q_n,P_n^*)-R(Q_0,P_0)+R(Q_0,P_{0,n}^*)+o_P(n^{-1/2})\\
&=& P_0 D^*_{Q_n}(P_n^*)-R_{2Q_n}(P_n^*,P_0)-P_0D^*_{Q_0}(P_{0,n}^*)+R_{2Q_0}(P_{0,n}^*,P_0)
+o_P(n^{-1/2})\\
&=& -(P_n-P_0)\{D^*_{Q_n}(P_n^*)-D^*_{Q_0}(P_{0,n}^*)\}\\
&&-\{R_{2Q_n}(P_n^*,P_0)-R_{2Q_0}(P_{0,n}^*,P_0)\}+o_P(n^{-1/2}). \end{eqnarray*}
At the last equality we used that $P_n D^*_{Q_n}(P_n^*)=P_nD^*_{Q_0}(P_{0,n}^*)=o_P(n^{-1/2})$. 
By the Donsker class condition, we have that the leading empirical process term is $O_P(n^{-1/2})$, and, by the third assumption, this shows that $d_0(Q_n,Q_0)=O_P(n^{-1/2})$. 
By our second assumption (where we can now use that $d_0(Q_n,Q_0)=O_P(n^{-1/2})$), and the asymptotic equicontinuity of the empirical  process indexed by a Donsker class, this yields that the leading empirical process term is $o_P(n^{-1/2})$.
This completes the proof. $\Box$

Regarding the first assumption we make the following remarks. 
Firstly, we note that if $P\rightarrow R(Q,P)$ would be a compactly differentiable functional so that $R(Q,P_n)$ is asymptotically linear estimator of $R(Q,P_0)$ for all $Q$, then the above lemma can be generalized to allow setting $P_n^*$ and $P_{0,n}^*$ equal to the empirical probability measure $P_n$. In that case, $R(Q_0,P_{0,n}^*)-R(Q_0,P_n^*)=0$, obviously.
In the case that $P\rightarrow R(Q,P)$ is non-smooth, the one  could set $P_n^*$ and $P_{0,n}^*$ equal to an undersmoothed HAL-MLE, so that again $P_n^*=P_{0,n}^*$ and thereby the first assumption holds. 
Finally, in the case that $P_n^*$ and $P_{0,n}^*$ are TMLEs targeting $R(Q_n,P_0)$ and $R(Q_0,P_0)$, respectively, one can still show that $R(Q_0,P_{0,n}^*)-R(Q_0,P_n^*)$ is a second order term so that assuming that it is $o_P(n^{-1/2})$ is a reasonable assumption.

\subsection{Efficiency theorem for the general HAL-MLE that allows for non-linear risk functions}
Consider again the general HAL-MLE w.r.t. a potentially non-linear risk function $P\rightarrow R(Q,P)$. We want to investigate conditions under which $\Psi(Q_n)$ is also asymptotically efficient for $\Psi(Q_0)$, thereby generalizing our Theorem \ref{thfinalspline} for the loss-based HAL-MLE.
Let $L_0(Q)(O)=R(Q,P_0)+D_Q(P_0)(O)$ and note that $L_0(Q)$ is a valid (when treating its dependence on $P_0$ as known)  loss function in the sense that $Q_0=\arg\min_{Q\in Q({\cal M})}P_0L_0(Q)$.

The following lemma provides sufficient conditions for
  $R(Q,P_n^*)=P_n L_0(Q)+r_n(Q)$ with a second order remainder $r_n(Q)$ that we can control uniformly in $Q\in Q({\cal M})$.
 \begin{lemma}\label{lemma1} 
 Recall $P_n^*=P_{n,Q}^*$ is such that $P_n D^*_Q(P_{n,Q}^*)=o_P(n^{-1/2})$ for all $Q$.
 Assume
 \begin{eqnarray*}
 \sup_{Q\in Q({\cal M})}P_n D^*_Q(P_{n,Q}^*)&=&o_P(n^{-1/2})\\
 \sup_{Q\in Q({\cal M})}\mid R_{Q,2}(P_n^*,P_0)\mid &=&o_P(n^{-1/2})\\
 \sup_{Q\in Q({\cal M})}\mid (P_n-P_0)\{D^*_Q(P_n^*)-D_Q(P_0)\}\mid&=& o_P(1).
 \end{eqnarray*}
Then,
\[
R(Q,P_n^*)=R(Q,P_0)+P_n D^*_Q(P_0)+r_n(Q),\]
where $\sup_{Q\in Q({\cal M})} \mid r_n(Q)\mid =o_P(n^{-1/2})$.
Here $r_n(Q)=(P_n-P_0)\{D^*_Q(P_n^*)-D^*_Q(P_0)\}+R_{2Q}(P_n^*,P_0)$.
\end{lemma}
{\bf Proof:} 
By $P_n^*=P_{n,Q}^*$ being a plug-in estimator satisfying $P_n D^*_Q(P_n^*)=o_P(n^{-1/2})$,  we have
\[
R(Q,P_n^*)-R(Q,P_0)=(P_n-P_0)D^*_Q(P_n^*)+R_{2Q}(P_n^*,P_0)+o_P(n^{-1/2}).\]
By assumption, $R_{2Q}(P_n^*,P_0)$ is $o_P(n^{-1/2})$ uniformly in all $Q$.
By assumption, the empirical process term equals $P_n D^*_Q(P_0)$ plus a negligible remainder, uniformly in $Q$.
Thus, we have
\[
R(Q,P_n^*)-R(Q,P_0)=P_n D^*_Q(P_0)+r_n(Q),\]
where $\sup_{Q\in Q({\cal M})}\mid r_n(Q)\mid =o_P(n^{-1/2})$.
$\Box$

This result suggests that minimizing $R(Q,P_n^*)$ is approximately the same as minimizing $P_n L_0(Q)$, where 
$L_0(Q)=R(Q,P_0)+D^*_Q(P_0)$ is an unknown loss function (i.e., a loss indexed by nuisance parameter). 
To further formalize this  we want to show that  the score equations $\frac{d}{d\epsilon}R(Q_{n,\epsilon}^h,P_{n,\epsilon}^*)$ at $\epsilon =0$ for $Q_n$ equal
the score equations $S_{h,0}(Q_n)\equiv \frac{d}{d\epsilon}P_n L_0(Q_{n,\epsilon}^h)$ up till an $o_P(n^{-1/2})$ approximation. With this result in hand, we can then simply apply Theorem \ref{thscoreequation} and Theorem  \ref{thfinalspline} with the loss function $L(Q)$ replaced by $L_0(Q)$, treating this loss function as given, and with the exact score equations $P_n S_h(Q_n)=0$  now replaced by $P_n S_{h,0}(Q_n)=o_P(n^{-1/2})$ for, uniformly in $h$ with $r(h,Q_n)=0$. This  provides us then with the conditions under which $\Psi(Q_n)$ is asymptotically efficient.  

Application of Lemma \ref{lemma1}  to $Q_{n,\epsilon}$ for a path $\{Q_{n,\epsilon}:\epsilon\}$ yields
\[
R(Q_{n,\epsilon},P_{n,\epsilon}^*)=R(Q_{n,\epsilon},P_0)+P_n D^*_{Q_{n,\epsilon}}(P_0)+r_n(Q_{n,\epsilon}).\]
It is reasonable to assume that $\left . \frac{d}{d\epsilon}r_n(Q_{n,\epsilon})\right |_{\epsilon =0}=o_P(n^{-1/2})$, because of two reasons. Firstly, $r_n(Q_{n,\epsilon})$ represents a second order remainder between $P_{n,\epsilon}^*-P_0$, where
$P_{n,\epsilon}^*$ is just a TMLE of $P_0$ targeting a particular target parameter and will thus converge just as fast as  the initial estimator $P_n^0$ used in the TMLE.
Secondly,  $r_n(Q_{n,\epsilon})$ is a second order remainder indexed by $Q_{n,\epsilon}$, so taking a derivative w.r.t. $\epsilon$ does not make this remainder worse.
This yields the following lemma.
\begin{lemma}\label{lemma2}
Recall definition of $r_n(Q)$ of Lemma \ref{lemma1}.
Define \[r_n(Q,h)\equiv\left .  \frac{d}{d\epsilon} r_n(Q_{\epsilon}^h)\right |_{\epsilon =0}.\]
Assume that for some upper bound $C<\infty$,
\[
\sup_{Q\in Q({\cal M}),\pl h\pl_{\infty}<C} \mid r_n(Q,h)\mid =o_P(n^{-1/2}).\]
Then, \[
\begin{array}{l}
\left. \frac{d}{d\epsilon}R(Q_{n,\epsilon},P_{n,\epsilon}^*) \right |_{\epsilon =0}=\left. \frac{d}{d\epsilon}R(Q_{n,\epsilon},P_0)\right |_{\epsilon =0}
+\left. P_n \frac{d}{d\epsilon}D_{Q_{n,\epsilon}}(P_0) \right |_{\epsilon=0}+o_P(n^{-1/2}).\end{array}
\]
\end{lemma}
Of course, this yields the following trivial corollary that is relevant for us.
\begin{corollary}
Assume the conditions of previous lemma.
Suppose that $\left . \frac{d}{d\epsilon}R(Q_{n,\epsilon}^h,P_{n,\epsilon}^*)\right |_{\epsilon =0}=0$ for a set  of paths indeed by $h\in {\cal H}$.
Then, 
\[
\left . P_n \frac{d}{d\epsilon} L_0(Q_{n,\epsilon})\right |_{\epsilon =0}=o_P(n^{-1/2})\mbox{
uniformly in all $h\in {\cal H}$.} \]
\end{corollary}
So we have now established conditions under which the general HAL-MLE is equivalent with a loss based HAL-MLE using loss function $L_0(Q)$, and that it also solves the same score equations of an HAL-MLE defined by minimizing $P_n L_0(Q)$ over $Q$.  
This allows us now to apply Theorem \ref{thscoreequation} and Theorem \ref{thfinalspline} with this choice of loss function.
This results in our next Theorem \ref{theffgenhalmle}.

Theorem \ref{theffgenhalmle} relies on the following definitions, analogue to Theorem \ref{thfinalspline}. \ \nl
{\bf Definitions:}
\begin{itemize} 
\item
Recall we can represent $Q_n=\arg\min_{Q\in {\cal Q}(C_n)}P_n R(Q,P_{n,Q}^*)$ as follows:
\[
Q_n(x)=I(Q_n)(x)+\sum_{\bar{s}(m)}\int_{(0_{s_m},x_{s_m}] } \phi_{\bar{s}(m),x_s}(u_{s_m})dQ^m_{n,\bar{s}(m)}(u_{s_m}),\]
where
\[
\begin{array}{l}
I(Q_n)(x)=Q_n(0)+\sum_{j=0}^{m-1}\sum_{\bar{s}(j)}Q_{n,\bar{s}(j)}^{j+1}(0_{s_j})\phi_{\bar{s}(j),\emptyset,x_s}(0_s).
\end{array}
\]
\item
 Consider the family of paths $\{Q_{n,\epsilon}^h:\epsilon\in (-\delta,\delta)\}$ through $Q_n$ at $\epsilon =0$ for arbitrarily small $\delta>0$, indexed by any uniformly bounded $h$, defined by
\begin{equation}\label{familypaths}
Q_{n,\epsilon}^h(x)=I(Q_{n,\epsilon}^h)(x)+\sum_{\bar{s}(m)}\int_{(0_{s_m},x_{s_m}] } \phi_{\bar{s}(m),x_s}(u_{s_m})(1+\epsilon h(\bar{s}(m),u_{s_m}))dQ^m_{n,\bar{s}(m)}(u_{s_m}),\end{equation}
where
\begin{eqnarray*}
I(Q_{n,\epsilon}^h)(x)&=&(1+\epsilon h(0))Q_n(0)
+\sum_{j=0}^{m-1}\sum_{\bar{s}(j)}\phi_{\bar{s}(j),\emptyset,x_s}(0_s)
(1+\epsilon h(\bar{s}(j),0_{s_j}))Q_{n,\bar{s}(j)}^{j+1}(0_{s_j}).
\end{eqnarray*}
\item  Let 
\begin{eqnarray*}
r(h,Q_n)&\equiv& I(h,Q_n)+\sum_{\bar{s}(m)}\int_{(0_{s_m},\tau_{s_m}] }h(\bar{s}(m),u_{s_m})\mid dQ^m_{n,\bar{s}(m)}(u_{s_m})\mid,\end{eqnarray*}
where
\[
\begin{array}{l}
I(h,Q_n)=h(0))\mid Q_n(0)\mid +\sum_{j=0}^{m-1}\sum_{\bar{s}(j)} h(\bar{s}(j),0_{s_j})  \mid Q_{n,\bar{s}(j)}^{j+1}(0_{s_j})\mid .
\end{array}
\]
\item For any uniformly bounded $h$ with $r(h,Q_n)=0$ we have that for a small enough $\delta>0$ $\{Q_{n,\epsilon}^h:\epsilon\in (-\delta,\delta)\}\subset {\cal Q}^m(C_n)$.\nl

\item Let $S_{0,h}(Q_n)=\left . \frac{d}{d\epsilon}L_0(Q_{n,\epsilon}^h)\right |_{\epsilon=0}$ be the $L_0(Q)$-score of this $h$-specific submodel
\item Consider the resulting set of scores 
\begin{equation}\label{scores}
{\cal S}_0(Q_n)=\{S_{0,h}(Q_n)=\frac{d}{dQ_n}L_0(Q_n)(f(h,Q_n)):\pl h\pl_{\infty}<\infty\},
\end{equation} where
\begin{eqnarray*}
f(h,Q_n)&\equiv &h(0)Q_n(0)+\sum_{s\subset \{1,\ldots,k\}}\int_{(0_s,x_s]} I((s,u_s)\in A) h(s,u_s) dQ_{n,s}(u_s).
\end{eqnarray*}
This is the set of scores generated by the above class of paths if we do not enforce constraint $r(h,Q_n)$.
\item 
Define 
\begin{eqnarray*}
r_n(Q)&\equiv& (P_n-P_0)\{D^*_Q(P_n^*)-D^*_Q(P_0)\}+R_{2Q}(P_n^*,P_0)\\
r_n(Q,h)&\equiv&\left .  \frac{d}{d\epsilon} r_n(Q_{\epsilon}^h)\right |_{\epsilon =0}.\end{eqnarray*}
Under assumption (\ref{critical}) below, the general HAL-MLE $Q_n=\arg\min_{Q\in Q({\cal M}),\pl Q\pl_v^*<C_n}R(Q,P_n^*)$ solves the score equations $P_n S_{0,h}(Q_n)=e_n(h)$,
up till an error $e_n(h)$ that is $o_P(n^{-1/2})$, uniformly in $h$, under constraint $r(h,Q_n)=0$: 
 $\sup_{\pl h\pl<\infty} \mid e_n(h)\mid =o_P(n^{-1/2})$. 
\item Let 
 \[
{\cal G}_n=\{G\in {\cal G}:D^*(Q_n,G)\in {\cal S}_0(Q_n)\} \]
be the set of $G$'s for which $D^*(Q_n,G)$ equals a score $S_{0,h}(Q_n)$ for some uniformly bounded $h$.
\item  Let $G_{0,n}\in {\cal G}_n$ be an approximation of $G_0$ so that $D^*(Q_n,G_{0,n})\in {\cal S}_0(Q_n)$.
Specifically, 
let $h^*(Q_n,G_0)$ so that $D^*(Q_n,G_{0n})=S_{0,h^*(Q_n,G_0)}(Q_n)$.
\item
Consider the representation (\ref{convenientrepresentation}) of $Q_n$:  $Q_n=\sum_{s,j\in {\cal J}_n(s)}\beta_n(s,j)\phi_{s,j}$ with $\beta_n$ a minimizer of the $R(\sum_{s,j\in {\cal J}_n(s)}\beta(s,j)\phi_{s,j},P_{\beta,n}^*)$ under constraint that its $L_1$-norm is bounded by $C_n$.
\end{itemize}

 \begin{theorem}\label{theffgenhalmle}
 Consider the above defined generalized HAL-MLE $Q_n=\arg\min_{Q\in {\cal Q}(C_n)}R(Q,P_n^*)$ for some $C_n\leq C^u<\infty$ with probability tending to 1. Consider also the above presented pathwise differentiable target parameter $\Psi:Q({\cal M})\rightarrow\openr^d$ with canonical gradient $D^*(Q(P),G(P))$ at $P\in {\cal M}$ and exact second order remainder $R_{20}((Q,G),(Q_0,G_0))=\Psi(Q)-\Psi(Q_0)+P_0D^*(Q,G)$.
\nl
{\bf Assumptions:}
\begin{itemize}
\item \begin{eqnarray*}
 \sup_{Q\in Q({\cal M})}\mid R_{Q,2}(P_n^*,P_0)\mid =o_P(n^{-1/2})\\
 \sup_{Q\in Q({\cal M})}\mid (P_n-P_0)\{D^*_Q(P_n^*)-D_Q(P_0)\}\mid&=& o_P(1).
 \end{eqnarray*}
\item For some upper bound $C<\infty$,
\begin{equation}\label{critical}
\sup_{Q\in Q({\cal M}),\pl h\pl_{\infty}<C} \mid r_n(Q,h)\mid =o_P(n^{-1/2}).\end{equation}

\item $P_0\{D^*_{Q_0}(P_n^*)-D^*_{Q_n}(P_n^*)\}^2$ can be bounded  by a constant $K_n$ times $d_0(Q_n,Q_0)$ with $K_n /n^{1/2}\rightarrow_p 0$.

\item $R_{20}(Q_n,G_{0n},Q_0,G_0)=o_P(n^{-1/2})$; $\sup_{Q\in Q({\cal M})}\mid R_{Q,2}(P_{n,Q}^*,P_0)\mid
=o_P(n^{-1/2})$;
$P_0\{D^*(Q_n,G_{0n})-D^*(Q_0,G_0)\}^2\rightarrow_p 0$;
$P_0\{D^*_{Q_0}(P_n^*)-D^*_{Q_n}(P_n^*)\}^2\rightarrow_p 0$.
Here we can use that $d_0(Q_n,Q_0)=R(Q_n,P_0)-R(Q_0,P_0)=o_P(n^{-1/2})$.
\item  $\{D_Q(P):P\in {\cal M},Q\in Q({\cal M})\}$ and $\{D^*(Q,G):Q\in {\cal Q}(C^u),G\in {\cal G}\}$ are contained in the class of $k_1$-variate cadlag functions on a cube $[0,\tau_o ]\subset\openr^{k_1}$ in a Euclidean space with a universal bound on the sectional variation norm.
\item $\lim\sup_n \pl h^*(Q_n,G_0)\pl_{\infty}<\infty$.
\item
\begin{equation}\label{assumptiona1}
\min_{s,j\in {\cal J}_n(s),\beta_n(s,j)\not =0}\pl P_n \frac{d}{dQ_n}L_0(Q_n)(\phi_{s,j})\pl =o_P(n^{-1/2}).\end{equation} 

\end{itemize}
Then, $d_0(Q_n,Q_0)=o_P(n^{-1/2})$; $P_n D^*(Q_n,G_0)=o_P(n^{-1/2})$, 
and $\Psi(Q_n)$ is an asymptotically efficient estimator of $\Psi(Q_0)$.
\end{theorem}

\end{document}